\documentclass{siamart220329}

\usepackage{amsmath}
\usepackage{amssymb}
\usepackage{amstext}
\usepackage{amscd}
\usepackage{arydshln}
\usepackage{color}
\usepackage{graphicx}
\usepackage{longtable}
\usepackage{calc}
\usepackage{microtype}
\usepackage{xcolor}
\usepackage{comment}
\usepackage{subcaption}
\usepackage{soul}

\newtheorem{remark}[theorem]{{\sc Remark}}

\newtheorem{example}[theorem]{{\sc Example}}

\newcounter{parenumi}
  {\end{list}}

\def\F{\mathbb{F}}

\def\R{\mathbb{R}}
\def\C{\mathbb{C}}
\def\N{\mathbb{N}}

\newlength{\tablecolwidth}
\definecolor{brilliantrose}{rgb}{1.0, 0.33, 0.64}
\definecolor{myviolet}{rgb}{0.21, 0.0, 0.85}
\definecolor{amethyst}{rgb}{0.6, 0.4, 0.8}
\definecolor{carrotorange}{rgb}{0.93, 0.57, 0.13}

\numberwithin{equation}{section}

\newcommand{\m}[1]{\begin{pmatrix} #1 \end{pmatrix}}
\newcommand{\bone}{\boldsymbol{1}}

\def\R{\mathbb{R}}
\def\C{\mathbb{C}}
\def\bone{{\bf 1}}

\def\N{\mathbb{N}}
\def\F{\mathbb{F}}
\DeclareMathOperator{\dd}{\textit{dd}}
\definecolor{carrotorange}{rgb}{0.93, 0.57, 0.13}

\makeatletter
\renewcommand*\env@matrix[1][*\c@MaxMatrixCols c]{%
  \hskip -\arraycolsep
  \let\@ifnextchar\new@ifnextchar
  \array{#1}}
\makeatother

\title{Nonbacktracking on Time-evolving Networks at the Node-level}
\author{Ryan Wood}
\date{March 2023}

\begin{document}

\title{Efficient computation of \lowercase{$f$}-centralities and nonbacktracking centrality for temporal networks}

\author{Vanni Noferini\thanks{Department of Mathematics and Systems Analysis, Aalto University,  P.O. Box 11100, FI-00076, Finland (\email{vanni.noferini@aalto.fi}).}
\and
Spyridon Vrontos\thanks{School of Mathematics, Statistics and Actuarial Science, University of Essex,  Wivenhoe Park,  Colchester, CO4 3SQ, United Kingdom
(\email{svrontos@essex.ac.uk}).}
\and
Ryan Wood\thanks{Corresponding author. Department of Mathematics and Systems Analysis, Aalto University,  P.O. Box 11100, FI-00076, Finland (\email{ryan.wood@aalto.fi}).}
}

\maketitle

\begin{abstract}
We discuss efficient computation of $f$-centralities and nonbacktracking centralities for time-evolving networks with nonnegative weights. We present a node-level formula for its combinatorially exact computation which proves to be more tractable than previously existing formulae at edge-level for dense networks. Additionally, we investigate the impact of the addition of a final time frame to such a time-evolving network, analyzing its effect on the resulting nonbacktracking Katz centrality. Finally, we demonstrate by means of computational experiments that the node-level formula presented is much more efficient for dense networks than the previously known edge-level formula. As a tool for our goals, in an appendix of the paper, we develop a spectral theory of matrices whose elements are vectors. 
\end{abstract}

\begin{keywords}
temporal network, f-centrality, nonbacktracking walk, matrices whose entries are vectors, Artin ring
\end{keywords}


\section{Introduction}
Complex networks appear in many fields \cite{strogatz2001exploring} and are a focal point for researchers and practitioners alike across a wide range of disciplines. There exists a multitude of network models, each seeking to encapsulate the different types of relationships and dependencies between actors within these networks. One such class are the networks with a time dependency. These time-evolving, or \emph{temporal}, networks have been the subject of much research over the past decade. In practice, their time-dependent structure typically manifests as a changing number of nodes, edges, or weights of edges. One common approach \cite{grindrod2011communicability} to modelling time-dependent networks has been to quantize the time-period over which we consider the network into discretely sampled time periods (also called time stamps or time frames), which are then treated as regular, non-time-dependent networks (also called static networks). Such a representation is sometimes called a `graph sequence representation' \cite{holme2015modern}. 

\textit{Centrality measures} are one of the most fundamental tools to analyze networks. The purpose of a centrality measure is to quantify the importance of each node within the network, and thereby identify the nodes most central to the network \cite{boldi2014axioms, estrada2012structure}. In practice, these are functions which assign to each node a non-negative value indicative of its centrality, within the network. There exist numerous families of centrality measures; this variety itself is reflective of the fact that ``centrality'', which may often be equated with ``importance'', is in practice a context-sensitive concept with a variety of interpretations. One broad class of centrality measures are those based on the combinatorics of walks. These centrality measures assign to each node a centrality value based on the weighted sum of the walks on the network that depart from, or end at, a given node \cite{benzi2013ranking}. A pioneering example of these centrality measures is \emph{Katz centrality} \cite{katz1953new}, that counts all walks that depart from a node, where walks of length $k$ are weighted as $t^k$ and $t$ is a positive parameter. By changing the weight to $c_k t^k$, where $(c_k)_k$ are the nonegative Taylor coefficients of a function $f(z)=\sum_{k=0}^\infty c_k z^k$, one can generalize Katz centrality to \emph{$f$-centrality} \cite{bb20,benzi2013ranking}; an important special case is exponential centrality, $c_k=\frac{1}{k!}$. More recently, another such walk-based centrality measure was introduced, namely \textit{nonbacktracking Katz centrality} \cite{arrigo2018non,grindrod2018deformed}. Nonbacktracking Katz centrality is defined similarly as classical Katz, but discounting walks that are \emph{backtracking}, i.e., contain a subsequences of nodes of the form $\cdots aba \cdots$; the resulting centrality measure has been proven to possess several tangible benefits compared to more traditional walk-based centrality measures \cite{alon2007non,grindrod2018deformed,lin2019non,timar2021approximating}.

Extending the definition of the centrality measures described above to temporal networks is not always an easy task. For Katz centrality, this can be done in an efficient multiplicative manner that leverages the static network case \cite{grindrod2011communicability}, but generally this is not always possible for $f$-centralities and nonbacktracking centralities \cite{arrigo2022dynamic}. With the goal of filling this gap, in \cite{arrigo2022dynamic,arrigo2024weighted}, a fully general method for the computation of both $f$-centralities and  nonbacktracking Katz centrality for a time-evolving network was introduced. The technique described in \cite{arrigo2022dynamic,arrigo2024weighted} is an edge-level approach, as it involves computing the source, target, and weight matrices \cite{arrigo2024weighted, kempton2016, NQ24} associated with either the line graph or the nonbackgracking line graph (that is, the graph whose adjacency matrix is the Hashimoto matrix \cite{hashimoto1990}) of the original graph. For a network having $n$ nodes, such line graphs can have up to $O(n^2)$ nodes; and in a time-evolving network, this difference may be persistent through several time frames. Hence, the algorithms of \cite{arrigo2022dynamic,arrigo2024weighted} can be prohibitively slow when the involved graphs are large and edge-dense.

In the present work, we show that both classical $f$-centralities and nonbacktracking Katz centrality of a temporal network can be computed entirely at the node-level, by using a block matrix. This leads to an improved computational efficiency, particularly for temporal graphs that are dense in their number of edges \cite{noferiniwood2024}, i.e., have $O(n^2)$ edges where $n$ is the number of nodes. One such application, that especially motivated this paper, was the design of an algorithm to conduct financial mathematics research building on the work of \cite{arslan2024portfolio,pozzi2013spread}. These papers describe an investment optimization method that involves computing centrality measures of a graph representing correlations between stock returns. We aimed to improve the method by using temporal networks, which happen to be typically very dense in edges. More details are given in Section~\ref{sec:FinanceExperiment}.

The paper is structured as follows. Section \ref{sec:bg} recalls the necessary background. In Section \ref{sec:efficient} we obtain our main results, Theorem \ref{thm:fcent} and Theorem \ref{theorem:main}, that respectively describe efficient formulae to compute $f$-centralities and nonbacktracking Katz centrality on a temporal network; Theorem \ref{theorem:update} shows how to cheaply update nonbacktracking Katz centrality when a new timeframe is added to a temporal network. Section \ref{sec:experiments} presents numerical experiments that demonstrate the effectiveness of the new approaches. Finally, for the sake of the exposition, we postpone to Appendix \ref{sec:matvec} the (somewhat technical) development of a spectral theory for matrices whose entries lie in a finite-dimensional vector space over a field. This may be interesting per se, and it is crucial to obtain some of our main results in Sections \ref{sec:efficient} and \ref{sec:updates}.

\section{Background}\label{sec:bg}
\subsection{Time-evolving networks and centrality measures}\label{subsec:time-evolving}
We begin by recalling the standard definitions of a \textit{static}, i.e., not evolving in time, network. This is modelled by a weighted digraph (possibly with loops).
\begin{definition}[Weighted digraph]\label{def:network}
    A directed graph $G = (V,E)$ is an ordered pair of sets, $V$, the set of nodes and $E \subseteq V \times V$, the set of edges between these nodes. We say there exists an edge $v_i \rightarrow v_j$ from node $v_i$ to node $ v_j \in V$, if $(v_i, v_j) \in E$.  Furthermore, in the case of a weighted graph, we also associate a weight function $\omega: E \rightarrow \mathbb{R}_{> 0}$ to $G$, and call $\omega(v_i,v_j)$ the \textit{weight} of the edge $(v_i,v_j) \in E$.
\end{definition}

Throughout this paper we will assume all graphs to be finite, and to have as their nodes the set $V = [n] := \{1, 2, \dots, n\}$. We also note in passing that unweighted graphs may simply be regarded as a weighted network with weight function $\omega \equiv 
1$.

We now turn our attention to the time-evolving counterpart \cite{grindrod2011communicability} of Definition \ref{def:network}. A time-evolving network consisting of $N$ time frames will be represented by a collection of weighted digraphs that all share the same set of vertices, but whose edges may change from one time frame to another. Time-evolving graphs are defined formally as follows.
\begin{definition}[Time-evolving network]\label{def:timeevolving}
A finite time-evolving network $\mathcal{G}$ is a finite collection of (possibly weighted) digraphs, 
$ \mathcal{G} := (G_{[1]}, \dots, G_{[N]})$ such that each constituent network $G_{[i]} = (V, E_{[i]})$ has an identical set of nodes but with a possibly varying set of edges 
and weight function. This sequence is associated with a non-decreasing times $(t_1, \dots, t_N) \in \mathbb{R}^{N}$ such that when observed at time $t_1$, $\mathcal{G}$ is identical to that of $G_{[i]}$.
\end{definition}
It is also useful to define the following type of subnetwork which often arises when treating time-evolving networks with a graph sequence representation.
\begin{definition}[Time-evolving subnetwork]
    Let $\mathcal{G} = (G_{[1]}, \dots, G_{[N]})$ be a time-evolving network defined as above. We call the network $\mathcal{H}:= (G_{[r]}, \dots, G_{[s])}$ a temporal subnetwork of $\mathcal{G}$ if, for all $1\leq r \leq i \leq s \leq N$, $G_{[i]}$ appears in both $\mathcal{G}$ and $\mathcal{H}$.
\end{definition}
In order to enumerate (with weights) walks within our temporal setting, we must first define what we mean by a walk across a finite time-evolving network. Intuitively, the idea is that we can walk around the nodes of the underlying graph, travelling over the existing edges (according to their direction), and taking into account that the configuration of the edges can change with time while we are still completing our walk. 
\begin{definition}[Walk of length $k$]\label{def:Twalk2}
{\rm A walk of length $k$} across a temporal network is defined as an ordered sequence of $k+1$ nodes $i_1 i_2 \dots i_{k+1}$ such that for all $\ell = 2,\ldots,k$ it holds that $i_{\ell-1}\to i_{\ell} \in E_{[r]}$ and $i_\ell\to i_{\ell+1} \in E_{[s]}$ for some $1\leq\tau_1 \leq r \leq s \leq \tau_2\leq N$, in which case the walk is said to begin on time frame $\tau_1$ and end on time frame $\tau_2$. Furthermore, the weight of such a walk is said to be equal to $\Pi_{\ell = 1}^{k-1} \omega_{[\tau_\ell]} ((i_\ell, i_{\ell+1}))$, i.e. the product of the weights of the constituent edges $i_\ell \rightarrow i_{\ell+1} \in E_{[\tau_\ell]}$ using the appropriate weight function for each time frame, $\omega_{[\tau_\ell]}$.
\end{definition}

Given a time-evolving network, we may represent each time frame with an \emph{adjacency matrix}, which encodes both the adjacency relations of nodes (whence the name) and the weighting of the edges.
\begin{definition}[Adjacency matrices associated with a time-evolving network]\label{def:adj}
Given a finite time-evolving network $\mathcal{G} = (G_{[1]}, \dots, G_{[N]})$. We associate with each $G_{[i]}$ the adjacency matrix $A_{[i]} \in \mathbb{R}^{n \times n}$, where $n$ is the constant number of nodes in the time-evolving network, such that
\[({A_{[i]}})_{v,w} = \begin{cases}
    \omega_{[i]}(v,w) & \textit{If }v\rightarrow w \in E_{[i]};\\ 0 & \textit{Otherwise.}
\end{cases}\]
    where $\omega_{[i]}$ is the weight function associated with $G_{[i]}$ described in Definition~\ref{def:network}.
\end{definition}

Powers of the adjacency matrix provide the weighted enumeration of the walks that occur across the network in a given timeframe, with weights as described in Definition~\ref{def:Twalk2}. We omit the elementary proof of the well known Lemma \ref{lem:walkcountingprop}.

\begin{lemma}\label{lem:walkcountingprop}
    Let $A_{[\ell]}$ be the adjacency matrix associated with the $\ell$-th timeframe of a time-evolving network $\mathcal{G}$, as defined in Definition~\ref{def:adj}, then the $(i,j)$-th entry of $A_{[\ell]}^k$ is equal to the sum of the weights of all walks which occur at timeframe $\ell$, begin from node $i$, and end in node $j$.
\end{lemma}
For a static network ($N=1$), Lemma \ref{lem:walkcountingprop} was utilized by Leo Katz in his seminal paper \cite{katz1953new} to define the following centrality measure, known now as \emph{Katz centrality}.

\begin{definition}[Katz Centrality]\label{def:Katzoriginal}
    The Katz centrality vector for a static network with $n$ nodes is
    \begin{equation}\label{eq:openkatz}
        \boldsymbol{x}_\textit{Katz}(t) := (I + \sum_{k =1}^\infty t^k A^k)\bone 
    \end{equation}
    where $\bone \in \R^n$ is the vector of all ones and $0<t<1$ is an attenuating factor which down-weights longer walks. If $0 < t < 1/\rho(A)$, where $\rho(A)$ is the spectral radius of $A$, then (\ref{eq:openkatz}) may be computed as
$\boldsymbol{x}_\textit{Katz}(t) = (I-tA)^{-1}\bone.$
\end{definition}
    Katz centrality of node $i$ equals to a weighted sum of all walks across the (static) network that begin from the node $i$; each walk of length $k$ is scaled by its weight as in Definition \ref{def:Twalk2} times $t^k$. Katz centrality was generalized to time-evolving networks via the formation of the so-called \emph{dynamic communicability matrix}, $\mathcal{Q}(t)$ \cite{grindrod2011communicability}.
\begin{definition}[Katz centrality for time-evolving networks]\label{def:dyncom}
    Given a time-evolving network consisting of $N$ time frames,  $\mathcal{G} = (G_{[1]}, \dots, G_{[N]})$, and a parameter $0< t < 1/\max_i(\rho(A_{[i]}))$, then the time-evolving Katz centrality of the network is \[\boldsymbol{x}_\textit{Katz}(t) = \mathcal{Q}(t) \bone := (I-tA_{[1]})^{-1}\cdots \dots \cdot  (I-tA_{[N]})^{-1} \bone. \]
\end{definition}
\subsection{Nonbacktracking walks and nonbacktracking Katz}
In Definition~\ref{def:Twalk2}, a walk on a time-evolving network was defined simply as a sequence of adjacent nodes that respect the time-evolving nature of the network, without any further sequential restrictions placed on the nodes that feature in the walk. One such restriction that has arisen as a focal point of research over the past decade (see, e.g., \cite{alon2007non,arrigo2024weighted,grindrod2018deformed,NQ24}) is that of \emph{nonbacktracking}; a walk that fails to adhere to this restriction is said to be \emph{backtracking} and is defined as follows.

\begin{definition}[Backtracking and nonbacktracking walk]\label{def:backtracking}
    A walk of length $k$, $i_1 i_2, \dots, i_{k+1}$ as described in Definition~\ref{def:Twalk2}
 is said to be {\rm backtracking} if for some $1 \leq l \leq k-1$, we have $i_l = i_{l+2}$, i.e. the walk involves a sequence of successive edges of the form $i_l \rightarrow i_{l+1}\in E_{[\tau_{1}]}$ and $i_{l+1} \rightarrow i_{l+2}\in E_{[\tau_{2}]}$ for some $1 \leq \tau_1 \leq \tau_2 \leq N$. A walk that is not backtracking is said to be \emph{nonbacktracking}.  
 \end{definition}
The nonbacktracking generalization of Katz centrality, that is, the weighted enumeration of nonbacktracking walks, was first obtained in \cite{arrigo2024weighted}. We recall it in Theorem \ref{thm:Phi}, which is stated using the somewhat lighter notation in \cite{NQ24} rather then its (equivalent, but more cumbersome) expression originally appeared in \cite{arrigo2024weighted}.
\begin{theorem}[Nonbacktracking Katz centrality for a weighted, static network]\cite{arrigo2024weighted, NQ24}\label{thm:Phi}
Let $P_k(A)$ be the matrix whose $(i,j)$-th entry is equal the sum of the weights of nonbacktracking walks, as per Definition~\ref{def:Twalk2}, departing from node $i$ and ending at $j$ of length $k$. Then, for all values of $t$ such that the sum $\sum_{k \geq 1} t^kP_k(A)$ converges, we have 

\begin{equation}\label{eq:weightedstatic}
I + \sum_{k \geq 1}t^k P_k(A) = (I - t\tilde{A}(t) + t^2 \tilde{D}(t))^{-1}
\end{equation}
where  $\tilde{D}(t)$ is diagonal and
\begin{align*}
    (\tilde{A}(t))_{ij} &:= \begin{cases}
    \omega((i,j)) & \textit{if }(i,j) \in E \textit{ and }(j,i) \not \in E \\
    \dfrac{\omega((i,j))}{1 - t^2\omega((i,j))\omega((j,i))} & \textit{if }(i,j) \in E \textit{ and } (j,i) \in E \\
    0 & \textit{otherwise,}
\end{cases}
\\
(\tilde{D}(t))_{ii} &:= \begin{cases}
    \sum_{j \in \Gamma(i) }\dfrac{\omega((i,j))\omega((j,i))}{1-t^2\omega((i,j))\omega((j,i))} & \textit{if } \ \exists \ j : (i,j) \in E \textit{ and } (j,i) \in E \\
    0 & \textit{otherwise.}
\end{cases}
\end{align*}
\end{theorem}

\section{Efficient computation of $f$-centralities and nonbacktracking centralities on time-evolving networks}\label{sec:efficient}
We now turn our attention to walk-based centrality measures on time-evolving networks. Let us begin with classical $f$-centralities, for which general computational method exist \cite{arrigo2022dynamic}. However, generally -- with notable exceptions \cite[Theorem 3.1]{arrigo2022dynamic}, such as $f(t)=(1-t)^{-1}$ that corresponds to Katz centrality -- these methods involve generating the adjacency matrix of the line graph, and/or other edge-level ``source" and ``target" matrices associated with each time frame as described for example in \cite{arrigo2020beyond,arrigo2022dynamic,kempton2016,NQ24}. The situation is similar for nonbacktracking Katz centrality, where the same edge-level method can be applied \cite{arrigo2022dynamic,arrigo2024weighted} except that the adjacency matrices of the line graphs are replaced by the \textit{Hashimoto matrix} \cite{hashimoto1990}.
However, the computational cost of this approach can be prohibitively high in the case of networks with many edges or many time frames. While sparse networks are more common, especially in sociological contexts, dense networks are not unheard of \cite{noferiniwood2024}. For example, in financial mathematics, certain techniques for portfolio selection \cite{pozzi2013spread} involve forming a network from the correlation matrix of a large number of stocks, which will often be (close to) the complete graph. Motivated by these facts, we investigate an alternative method using a node-level adjacency matrix, that can be far more efficient while being mathematically equivalent to the previously known formulae (see \cite[Section 5]{arrigo2022dynamic} for classical $f$-centralities and Theorem \ref{thm:Phi} above for nonbacktracking Katz centrality).

\subsection{The time-evolving adjacency matrix}
In this section we begin by describing the basic objects that we will associate with a time-evolving network that will facilitate the enumeration of (nonbacktracking) walks with weights. The first of these is what could be called the time-evolving analogue of the adjacency matrix.

\begin{definition}[Time-evolving adjacency matrix]{\label{def:globalA}}
    Consider a finite time-evolving graph $\mathcal{G} =  (G_{[1]}, \dots, G_{[N]})$ with adjacency matrices $A_{[i]}$ as in Definition \ref{def:adj}. The node-level adjacency matrix of $\mathcal{G}$ is the block upper-triangular matrix
    \begin{equation}\label{eq:teae}\mathcal{A} = \m{A_{[1]} & A_{[2]} & \dots & A_{[N]} \\ 0 & A_{[2]} & \dots & A_{[N]} \\
    0 & 0 & \ddots &  A_{[N]}\\
    }.  \end{equation}
\end{definition}

The matrix \eqref{eq:teae} could be viewed as a special case of a more general framework of block matrices for temporal networks previously introduced in the literature; see, e.g., \cite{al2021block}. Here, we take a combinatorial viewpoint and we add the crucial observation that, by specializing to \eqref{eq:teae}, we retain the familiar walk-counting property of adjacency matrices associated with static graphs. 
Indeed, Lemma \ref{lem:walkcount} generalizes Lemma \ref{lem:walkcountingprop}.
\begin{lemma}[Walk enumeration property]\label{lem:walkcount}
Let $\mathcal{A}$ as in~\eqref{eq:teae} be the node-level adjacency matrix of a temporal network $\mathcal{G} = (G_{[1]}, \dots, G_{[N]})$ with $n$ nodes. For all $0 < k \in \mathbb{N}$, $1 \leq \tau_1 \leq \tau_2 \leq N$ and $1 \leq i,j \leq n$, the $((\tau_1-1) n + i,(\tau_2-1) n + j)$ entry of $\mathcal{A}^k$ is the sum of the weights of walks (as per Definition~\ref{def:Twalk2}) of  length $k$ which begin from node $i$ and terminate in node $j$, for which the first edge of the walk is in $E_{[\tau]}$, where $\tau_1 \leq \tau \leq \tau_2$, and for which the final edge of the walk is in $E_{[\tau_2]}$.
\end{lemma}
\begin{proof}
Let $\mathcal{G}$, $1 \leq \tau_1 \leq  \tau_2 \leq N$ and $1 \leq i, j \leq n$ be given. We first alleviate the notation by referring to the block-structure inherited from $\mathcal{A}$. In particular we say 
\[{\mathcal{A}^k}_{((\tau_1-1) n + i,(\tau_2-1) n + j)} =: ({\mathcal{A}^k}_{[\tau_1, \tau_2]})_{i,j}, \]
where ${\mathcal{A}^k}_{[\tau_1, \tau_2]}$ may be called the $(\tau_1,\tau_2)$-th block of $\mathcal{A}^k$. We now proceed by induction on $0 < k \in \mathbb{N}$. For the base case, we examine the walks of length $1$ beginning from node $i$ and terminating in node $j$, which satisfies the additional requirement that the first edge belongs to $E_{[\tau]}$ and the final edge belongs to $E_{[\tau_2]}$, where $\tau_1 \leq \tau \leq \tau_2$. Since the walk must be of length $1$, then $\tau = \tau_2$, and thus the unique such walk is $w:=i\xrightarrow[]{\tau_2}j \in E_{[\tau_2]}$. But the $(\tau_1, \tau_2)$-th block of $\mathcal{A}$ is equal to $A_{[\tau_2]}$, which by Definition~\ref{def:adj} has as its $(i,j)$-th entry $k = 1$ the weight of the walk $w$. Assume now that the statement holds for all walks of length $1,2,\dots, k$, with the aim to prove the statement for walks of length $k+1$. Examining the relevant entry of $\mathcal{A}^{k+1}$ yields

\begin{align*}
(\mathcal{A}^{k+1})_{(\tau_1-1) n + i,(\tau_2-1) n + j} =\left((\mathcal{A}^{k+1})_{[\tau_1, \tau_2]}\right)_{i,j} &= \left( \sum_{m = 1}^N \mathcal{A}_{[\tau_1,m]} (\mathcal{A}^k)_{[m,\tau_2]} \right)_{i,j}\\
&= \left(\sum_{m = \tau_1}^{\tau_2} A_{[m]} (\mathcal{A}^k)_{[m,\tau_2]}\right)_{i,j}.
\end{align*}
The final summand can be interpreted as extending all walks of length $k$, which start on or after the $m$-th timeframe and end on the timeframe $\tau_2 \geq m$, by a single edge taken from $E_{[m]}$. Therefore, by summing over $m$ which takes values from $\tau_1 $ to $\tau_2$, we obtain the sum of the weights of all the walks of length $k+1$ which begin at some time between $\tau_1$ and $\tau_2$ which terminate during the $\tau_2$-th timeframe. 
\end{proof}

\subsection{Classical $f$-centralities}
Following \cite{arrigo2022dynamic}, given a function $f(t)=\sum_{k=0}^\infty c_k t^k$, analytic at $t=0$ and such that $c_k \geq 0$, we define the $f$-centrality of node $i$ in a temporal network as the weighted sum of all walks (as in Definition \ref{def:Twalk2}), where the weight of a walk of length $k$ is $c_k t^k$. Using Lemma \ref{lem:walkcount}, we can easily compute the $f$-centrality for a time-evolving network, where $f$ is an analytic function with non-negative Taylor coefficients.

\begin{theorem}[$f$-centrality via $\mathcal{A}$] \label{thm:fcent} Let $\mathcal{G} = (G_{[1]}, \dots, G_{[N]})$ be a time-evolving network. Given an analytic function $f(z) = \sum_{i = 0}^\infty c_i z^i$ for which all coefficients $c_i \geq 0$ and radius of convergence $r$, the $f$-centrality of the $i$-th
node for the time-evolving subnetwork $\mathcal{G} = (G_{[\tau_1]} \dots, G_{[\tau_2]})$, $1 \leq \tau_1 \leq \tau_2 \leq N$ is given by
\[\sum_{j = n(\tau_1-1) + 1}^{\tau_2 n} [f(t \mathcal{A})]_{(\tau_1-1)n + i,j} \]
 where $t$ is the attenuation factor chosen from the interval $0 < t <r/\rho(\mathcal{A}) $.
\end{theorem}
\begin{proof}
Consider the series
    \begin{equation}\label{eq:theseries}
         f(t\mathcal{A}) = \sum_{k = 0}^\infty c_k t^k \mathcal{A}^k, 
    \end{equation}
    that by classical linear algebra results has radius of convergence $r [\rho(\mathcal{A})]^{-1}$. We have from Lemma~\ref{lem:walkcount} that $\mathcal{A}^{k}$ enumerates in its $[\tau_i, \tau_j]$-th block the weighted sum of walks of length $k$ that begin on some time frame $\tau$ where $\tau_i \leq \tau \leq \tau_j$ and end on the time frame $\tau_j$. Hence, \eqref{eq:theseries}
    has as its $[\tau_i, \tau_j]$-block the sum of weighted walks further weighted by $c_k t^k$ (where $k$ is the length of the respective walk being summed) that begin on some time frame $\tau$ where $\tau_i \leq \tau \leq \tau_j$ and end on the time frame $\tau_j$. By fixing the (scalar) row index to $(\tau_1-1)n+i$, and by summing over the (scalar) column index $j$ between $n(\tau_1-1)+1$ and $\tau_2 n$, we thus obtain the weighted sum of all walks which begin from node $i$ on some time frame $\tau \geq \tau_1$ but end at the latest during the $\tau_2$-th time frame. This is precisely the temporal $f$-centrality of the node $i$ in the temporal network $(G_{[\tau_1]} \dots, G_{[\tau_2]})$.
\end{proof}
By specializing Theorem \ref{thm:fcent} to $f(t)=\frac{1}{1-t}=\sum_{k=0}^\infty t^k$, we recover known results for dynamic Katz centrality. In particular, since summing by columns is equivalent to right-multiplying by $\bone$, we have the following Corollary.

\begin{corollary}\label{corollary:resolvent}
Let $\mathcal{G} := (G_{[1]}, \dots,  G_{[N]})$ be a time-evolving network with fixed nodes $V = [n]$. Further let $\mathcal{A}$ be the time-evolving adjacency matrix associated with $\mathcal{G}$. Then time-evolving Katz centrality of the node $i \in V$ within the time-evolving subnetwork $\mathcal{H} = (G_{[s]}, \dots, G_{[N]})$ is given by the following:
    \[((I - t \mathcal{A})^{-1} \boldsymbol{1})_{((s-1)n +i)}\]
    where $t$ is the attenuation factor chosen from the interval $0 < t <1/\rho(\mathcal{A}) $.
\end{corollary}
    The following example illustrates that the expression in Corollary \ref{corollary:resolvent} is equivalent to the familiar dynamic communicability matrix as per Definition~\ref{def:dyncom}.\begin{example}\label{example:Katz}
    Let $\mathcal{G} = (G_{[1]}, G_{[2]}, G_{[3]})$, then per Definition~\ref{def:globalA}, the time-evolving adjacency matrix associated with $\mathcal{G}$ is the following,
    \[\mathcal{A} = \m{A_{[1]} & A_{[2]} & A_{[3]} \\ 0 & A_{[2]} & A_{[3]} \\ 0 & 0 & A_{[3]}}\]
    If we then calculate the resolvent function $X := (I-t\mathcal{A})^{-1}$, then we obtain the following submatrix blocks:
    \begin{align*}
    X_{[1,1]} &= (I-tA_{[1]})^{-1} \\
    X_{[1,2]} &= (I-tA_{[1]})^{-1}tA_{[2]}(I-tA_{[2]})^{-1} =  (I-tA_{[1]})^{-1} ((I-tA_{[2]})^{-1} - I)\\
    X_{[1,3]} &=  (I-tA_{[1]})^{-1}(I-tA_{[2]})^{-1}tA_{[3]}(I-tA_{[3]})^{-1} = (I-tA_{[1]})^{-1}(I-tA_{[2]})^{-1}((I-tA_{[3]})^{-1} - I)
    \end{align*}
The sum of these then gives
\begin{equation*}
    X_{[1,1]} + X_{[1,2]} + X_{[1,3]} = (I-tA_{[1]})^{-1}(I-tA_{[2]})^{-1}(I-tA_{[3]})^{-1}=\mathcal{Q}
\end{equation*}
where $\mathcal{Q}$ is the dynamic communicability matrix (Definition~\ref{def:dyncom}).
\end{example}

\subsection{Nonbacktracking Katz centrality}\label{sec:NBTDynamic}

 In this subsection we examine how to extend Theorem~\ref{thm:fcent} to the case of \textit{nonbacktracking Katz centrality}. In other words, we will explain why the time-evolving adjacency matrix of Definition~\ref{def:globalA} may also be used to effiicently compute the nonbacktracking Katz centrality of weighted, possibly directed, temporal network.

Nonbacktracking Katz centrality is known to offer various advantages \cite{arrigo2018non, arrigo2018exponential, arrigo2020beyond,benzi2013ranking} over classical Katz centrality. However previously described methods for computing the nonbacktracking Katz centrality of a time-evolving network  \cite{arrigo2022dynamic, arrigo2024weighted} which can be computationally expensive when the network is not sparse. This motivates a new approach, which is largely inspired by the analysis of static weighted networks in \cite{arrigo2024weighted}. Namely, we  establish a recurrence relation which will we subsequently solve for a nonbacktracking generating function. The underlying mathematical tool is a spectral theory of $R^{N \times N}$ where $R$ is the commutative ring of square matrices over $\R \subset \C$, equipped with addition and elementwise multiplication; see Appendix \ref{sec:matvec} for a detailed treatment. We will use the notation developed in Appendix \ref{sec:matvec}, and in particular matrix multiplication over $R^{N \times N}$ will be denoted by $A \ast B$. As observed in Section \ref{sec:dee}, even if one can also interpret $A,B$ as elements of $\R^{nN \times nN}$, generally $A \ast B \neq AB$ where $AB$ is traditional matrix multiplication. For our goals in the present section, it is also convenient to define two additional operations on $R^{N \times N}$:
    \begin{definition}
     \[If \ {A}^{\ast T} := \m{A_{11} & A_{12} \\ A_{21} & A_{22}},\ then \ {A}^{\ast T} := \m{A_{11}^T & A_{12}^T \\ A_{21}^T & A_{22}^T}\]
    where $A^T$ is the regular matrix transpose of $A$,
   and
    \[\dd^{\ast}(A) := \m{\dd(A_{11}) & \dd(A_{12}) \\ \dd(A_{21}) & \dd(A_{22})} \]
    where $\dd(A)$ is the diagonal matrix whose diagonal equals that of $A$.
    \end{definition}

In order to remove the walks that backtrack, we first relate the powers (in $R^{N \times N}$) of $\mathcal{A} \ast \mathcal{A}^{\ast T}$ to the number of walks that use only specific nodes.
\begin{lemma}\label{lemma:astcounting}
Let $\mathcal{A}$ be as in Definition \ref{def:globalA}. Then,
    $(\mathcal{A} \ast \mathcal{A}^{\ast T})^{\ast k}$ counts in its $((\tau_1-1)n + i,(\tau_2-1)n + j)$-entry the number of walks of length $2k$ using only nodes $i$ and $j$ which end on time frame $\tau_2$ but begin on some time frame $\tau_1 \leq \tau \leq \tau_2$.
\end{lemma}
\begin{proof}
By definition,
\[\left( \mathcal{A}\ast \mathcal{A}^{\ast T}\right)_{\tau_1, \tau_2}= \sum_{r = 1}^{N} \mathcal{A}_{\tau_1, r} \circ {\mathcal{A}^{\ast T}}_{ r, \tau_2} =\sum_{r = \tau_1}^{\tau_2} A^{[ r]} \circ {A^{[\tau_2]}}^T \]
where the second equality follows from the upper-triangular block-structure of $\mathcal{A}$. The $(i,j)$-th entry of this sum counts the number of walks of the form
\[i \xrightarrow{\tau_1 \leq \boldsymbol{r} \leq \tau_2} j \xrightarrow{ \boldsymbol{\tau_2}} i\]
where the symbols in bold indicated the timeframe to which each edge belongs.

 Having proven that the result holds for this base case of $k = 1$, assume now that the result holds for all $r < k$, so that 
$(\mathcal{A} \ast \mathcal{A}^{\ast T})^{\ast (k-1)}_{\tau_1,\tau_2}$ counts in its $(i,j)$-th entry the number of walks of length $2(k-1)$ that alternate from node $i$ to node $j$, and that begin at the earliest on time frame $\tau_1$ but end on time frame $\tau_2$. Now, express the $(\tau_1,\tau_2)$ block of $(\mathcal{A} \ast \mathcal{A}^{\ast T})^{\ast k}$ as the following sum:

\begin{equation}
  \label{eq:followingsum}
 \left[(\mathcal{A} \ast \mathcal{A}^{\ast T})^{\ast k}\right]_{\tau_1,\tau_2}= \sum_{s = \tau_1}^{\tau_2}\left( \sum_{r = \tau_1}^{s} A^{[r]}\circ {A^{[s]}}^T\right)\circ \left((\mathcal{A} \ast \mathcal{A}^{\ast T})^{\ast (k-1)}\right)_{s,\tau_2}.
\end{equation}

The first summand $ \sum_{r = \tau_1}^{s} A^{[r]}\circ {A^{[s]}}^T$ in \eqref{eq:followingsum} counts in its $(i,j)$-th entry the number of walks of the form $i \xrightarrow{r} j\xrightarrow{ s} i $, where  $\tau_1 \leq r \leq s \leq \tau_2$. Thus, overall \eqref{eq:followingsum} counts the number of walks of the form
 \[i \xrightarrow{\boldsymbol{r} \geq \tau_1  } j \xrightarrow{\boldsymbol{s}\geq r } i \xrightarrow{\boldsymbol{t} \geq s} j \xrightarrow{} \dots \xrightarrow{\boldsymbol{\tau_2}} i \]
The result is now proven, it remains only to note that the $(i,j)$-th entry of the $(\tau_1, \tau_2)$-th block is the $((\tau_1-1)n + i,(\tau_2-1)n + j)$-th entry of $(\mathcal{A} \ast \mathcal{A}^{\ast T})^{\ast k}$.
\end{proof}
Lemma~\ref{lemma:astcounting} is a fundamental result, and we shall use it to derive a recurrence relation which, in turn, will allow us to enumerate the nonbacktracking walks. It allows us to correctly weight walks which span multiple time frames. 
We now define the matrices $P_k(\mathcal{A}) \in R^{N \times N}$ as the matrices whose $[\tau_i, \tau_j]$-th entry\footnote{To clarify, this is an entry when $P_k(\mathcal{A})$ is viewed as an element of $R^{N \times N}$. Alternatively,  it is an $n \times n$ block entry when the same matrix is viewed as an element of $\mathbb{R}^{nN \times nN}$.} is the matrix which counts in its $(i,j)$-th entry the total number of nonbacktracking walks of length $k$ which begin at the earliest on the time frame $\tau_i$ from node $i$ but end on the time frame $\tau_j$ in node $j$.
We further define $P_0(\mathcal{A}) = I$. With these definitions, we have the following recurrence relation.
\begin{theorem}[Recurrence Relation]\label{theorem:NBTrecurrence}
    For all $k \geq 1$,
\[P_k(\mathcal{A}) = \sum_{\substack{\ell=2h+1 \ \mathrm{odd} \\ 1 \leq \ell \leq k}} ((\mathcal{A} \ast \mathcal{A}^{\ast T})^{\ast h} \ast \mathcal{A}) P_{k-\ell}(\mathcal{A}) - \sum_{\substack{\ell=2h \ \mathrm{even} \\ 2 \leq \ell \leq k}} \dd^{\ast}(\mathcal{A} ((\mathcal{A} \ast \mathcal{A}^{\ast T})^{\ast h} \ast \mathcal{A})  ) P_{k-\ell}(\mathcal{A}).  \]
where $\mathcal{A}^{\ast h} = \mathcal{A} \underbrace{\ast \dots \ast}_{\textit{h times}} \mathcal{A}$, and $\mathcal{A}^{\ast 0 }$ is understood to be $E$, the identity element of the ring $R^{N\times N}$ described in Appendix~\ref{sec:dee}.
\end{theorem}

\begin{proof}
    We proceed by induction on $k$. The base case is trivial as $P_1(\mathcal{A}) = \mathcal{A} I = \mathcal{A}$, which indeed counts the number of nonbacktracking walks as specified. Now assume the statement holds for $P_\ell(\mathcal{A})$ whenever $\ell \leq k$, and we will prove that it also holds for $P_{k+1}(\mathcal{A})$. Consider a given nonbacktracking walk of length $k$ beginning from node $i$ and ending on node $j$ that ends at some timeframe $\tau$, $i \rightarrow \dots \xrightarrow{\boldsymbol{\tau}} j$. Such a walk can be extended on the left by left-multiplying by $\mathcal{A}$, the resulting walks are of the form
    \[i \xrightarrow{  \leq \tau} a \rightarrow b \dots \xrightarrow{\boldsymbol{\tau}} j\]
    and are counted with weights by $\mathcal{A}P_{k-1}(\mathcal{A})$. This walk is now possibly backtracking in its first two steps, i.e. it may be of the form
    \[i \xrightarrow{\boldsymbol{\tau_1} \leq \tau} a \xrightarrow{\tau_1 \leq \boldsymbol{\tau_2}\leq \tau} i \dots \xrightarrow{\boldsymbol{\tau}} j\]
    such walks\footnote{Note here that the notation $i \xrightarrow{\boldsymbol{\tau_1} \leq \tau} a$ is taken to mean that the edge $i\rightarrow a \in E^{[\tau_1]}$, where $\tau_1 \leq \tau$, i.e. the bold symbol is used to indicate the time frame to which the edge belongs.} are counted by $\dd^{\ast}( \mathcal{A} \mathcal{A})P_{k-2}(\mathcal{A})$. However $\dd^{\ast}( \mathcal{A} \mathcal{A})P_{k-2}(\mathcal{A})$ additionally counts walks of the form
    \[i \xrightarrow{\boldsymbol{\tau_1} \leq \tau} a \xrightarrow{\tau_1 \leq \boldsymbol{\tau_2} \leq \tau} i \xrightarrow{\tau_2 \leq \boldsymbol{\tau_3} \leq \tau} a \rightarrow \dots \xrightarrow{\boldsymbol{\tau}} j\]
    which were not included in $\mathcal{A}P_{k-1}(\mathcal{A})$; these walks are counted by $(\mathcal{A} \ast \mathcal{A}^{\ast T} \ast \mathcal{A})P_{k-3}(\mathcal{A})$. However this quantity again includes backtracking walks that haven't been counted by previous terms of the form:
    \[i \xrightarrow{\boldsymbol{\tau_1} \leq \tau} a \xrightarrow{\tau_1 \leq \boldsymbol{\tau_2} \leq \tau} i \xrightarrow{\tau_2 \leq \boldsymbol{\tau_3} \leq \tau} a \rightarrow  \xrightarrow{\tau_3 \leq \boldsymbol{\tau_4} \leq \tau} i  \rightarrow \dots \xrightarrow{\boldsymbol{\tau}} j\]
    These walks are counted by $\dd^{\ast} \left( \mathcal{A} (\mathcal{A} \ast \mathcal{A}^{\ast T}\ast  \mathcal{A}) \right)P_{k-4}(\mathcal{A})$, however this again overcounts in a similar manner as before. Continuing this procedure throughout the entire walk, we obtain the stated recurrence relation.
\end{proof}
In order to compute the nonbacktracking Katz centrality, we need to compute the generating function for the $P_k(\mathcal{A})$ matrices, i.e. $\Psi(\mathcal{A},t):=\sum_{k = 0}^\infty t^k P_k(\mathcal{A})$. This is achieved by first defining the generating function of the coefficients that appear in the recurrence relation of Theorem~\ref{theorem:NBTrecurrence}$, \Phi(\mathcal{A},t) = \sum_{h = 0}^\infty t^h C_h$, where
 \[C_h = \begin{cases} I & \textit{If } h = 0,\\ \dd^{\ast}(\mathcal{A} ((\mathcal{A} \ast \mathcal{A}^{\ast T})^{\ast h} \ast \mathcal{A})  ) & \textit{for } 0 \neq h \textit{ even,} \\-(\mathcal{A} \ast \mathcal{A}^{\ast T})^{\ast h} \ast \mathcal{A} & \textit{for }h \textit{ odd.}
    \end{cases}\]

Propositions \ref{prop:s} and \ref{prop:even} take care separately of the odd and even terms in $\Phi(\mathcal{A},t)$. To state them, it is convenient to define
\begin{equation}\label{eq:t0}
t_0(\mathcal{A}):=\left(\max_{i,j,s} \left\{(A_{[s]})_{ij} (A_{[s]})_{ji}  \right\} \right)^{-1/2}.
\end{equation}
\begin{proposition}\label{prop:s}
    For all $0 < t < t_0(\mathcal{A})$ as in \eqref{eq:t0} 
        we have
    \[\sum_{h = 0}^\infty t^{2h+1}(\mathcal{A} \ast \mathcal{A}^{\ast T})^{\ast h} \ast \mathcal{A} =   t \mathcal{A} \ast (E - t^2 \mathcal{A}^{\ast T} \ast \mathcal{A})^{\ast (-1)}.\]
\end{proposition}
\begin{proof}
   The sum can be rewritten as
\[ S=t \mathcal{A} \ast \left( \sum_{h=0}^\infty (t^2)^h (\mathcal{A}^{\ast T} \ast \mathcal{A})^{\ast h} \right).  \]
Setting $\mathcal{B}:=\mathcal{A}^{\ast T} \ast \mathcal{A}$, then by Theorem \ref{thm:68} the series within brackets has radius of convergence (using the notation introduced in Section \ref{sec:ee} of the appendix) \[ \sqrt{\frac{1}{\rho_R(\mathcal{B)}}} = \sqrt{\frac{1}{\max_{i,j}\rho( [\mathcal{B}]_{ij}  )}}.\] 
Hence, using also Corollary \ref{cor1} and the structure of $\mathcal{A}$, \[ ([\mathcal{B}]_{ij})_{rs}=(\mathcal{B}_{rs})_{ij} = ([\mathcal{A}^{\ast T}]_{ij}[\mathcal{A}]_{ij}  )_{rs} = \sum_{k=1}^n ([\mathcal{A}^{\ast T}]_{ij})_{rk} ([\mathcal{A}]_{ij}  )_{ks} = \sum_{k=r}^s  (A_{[k]}^T)_{ij} (A_{[s]})_{ij}.\]
This shows in particular that $[\mathcal{B}]_{ij}$ is upper triangular and that its $(s,s)$ diagonal element is $(A_{[s]})_{ij}(A_{[s]})_{ji}$, thus yielding the first part of the statement. On the other hand,
\[ S  = t \mathcal{A} + t^3 \mathcal{A} \ast \mathcal{A}^{\ast T} \ast \mathcal{A} + t^5\mathcal{A} \ast \mathcal{A}^{\ast T} \ast \mathcal{A} \ast \mathcal{A}^{\ast T} \ast \mathcal{A} + \dots     \]
and hence  $S \ast t^2(\mathcal{A}^{\ast T} \ast \mathcal{A} ) = S-t\mathcal{A}\iff t \mathcal{A} = S \ast (E - t^2 \mathcal{A}^{\ast T}\ast \mathcal{A}).$
\end{proof}

\begin{proposition}\label{prop:even}
    For all $0 < t < t_0(\mathcal{A})$ as in \eqref{eq:t0}    we have
    \[  \sum_{h=0}^\infty     t^{2h+2}\dd^\ast(\mathcal{A} ((\mathcal{A} \ast \mathcal{A}^{\ast T})^{\ast h} \ast \mathcal{A})  )  = \dd^{\ast}\left( t\mathcal{A}(t \mathcal{A} \ast (E - t^2 \mathcal{A}^{\ast T} \ast \mathcal{A})^{\ast (-1)}  )\right).    \]
\end{proposition}
\begin{proof}
Using also Proposition \ref{prop:s},   we have that the assumption on $t$ implies that \begin{align*}
        \sum_{h = 0}^\infty t^{2h+2}\dd^\ast(\mathcal{A} ((\mathcal{A} \ast \mathcal{A}^{\ast T})^{\ast h} \ast \mathcal{A})  )  &= \dd^{\ast}\left[ t\mathcal{A} \left(  \sum_{h = 0}^\infty t^{2h+1}(\mathcal{A} \ast \mathcal{A}^{\ast T})^{\ast h} \ast \mathcal{A} \right)\right] \\
        &=\dd^{\ast}\left( t\mathcal{A}(t \mathcal{A} \ast (E - t^2 \mathcal{A}^{\ast T} \ast \mathcal{A})^{\ast (-1)}  )\right). 
    \end{align*}
\end{proof}

\begin{theorem}\label{theorem:main}
    Let $\Psi(\mathcal{A},t)=\sum_{k = 0}^\infty t^k P_k(\mathcal{A})$. Then, for all $0 < t < t_0(\mathcal{A})$ as in \eqref{eq:t0}
    it holds
    \[\Psi(\mathcal{A},t)=
    \left[ I - t \mathcal{A} \ast (E - t^2 \mathcal{A}^{\ast T} \ast \mathcal{A})^{\ast (-1)}  + \dd^\ast\left( t\mathcal{A}(t \mathcal{A} \ast (E - t^2 \mathcal{A}^{\ast T} \ast \mathcal{A})^{\ast (-1)}  )\right) \right]^{-1} .\]
\end{theorem}
\begin{proof}
    By \cite[Proposition 3.2]{arrigo2024weighted}, $ I = \Phi(\mathcal{A},t) \Psi(\mathcal{A},t) = (R-S+I) \Psi(\mathcal{A},t),$
    where $R$ (resp. $S$) denote the infinite sum appearing in the statement of Proposition \ref{prop:even} (resp. Proposition \ref{prop:s}). The statement follows immediately.
\end{proof}

From Theorem~\ref{theorem:main}, we can compute the temporal nonbacktracking Katz centrality.

\begin{corollary}\label{corollary:howtocompute}
    Let $\mathcal{G} = (G_{[1]}, \dots, G_{[N]})$ be a time-evolving network with $n$ nodes as per Definition~\ref{def:timeevolving} and $\mathcal{A}$ its associated time-evolving adjacency matrix. For all $0 < t < t_0(\mathcal{A})$ as in \eqref{eq:t0}, the nonbacktracking Katz centrality of $\mathcal{G}$ is given by the top $n$ entries of the vector $\Psi(\mathcal{A},t)\bone$ where $\Psi(\mathcal{A},t)$ is as in Theorem \ref{theorem:main}. 
\end{corollary}
\begin{remark}
    More generally, the nonbacktracking Katz centrality of any time-evolving subnetwork $\mathcal{H}$ of $\mathcal{G}$ can be computed via the application of Corollary~\ref{corollary:howtocompute} to the relevant submatrix of $\Psi(\mathcal{A})$ and multiplying this by the vector of all ones $\bone$. For example, when $\mathcal{H} = (G_{[r]}, \dots, G_{[N]})$, the centrality of the $i$-th node is given by the $(r-1)n + i$-th entry of $\Psi(\mathcal{A},t)\bone$
\end{remark}
   For a single time frame, i.e., $\mathcal{G} = (G_{[1]})$, we see that  $\mathcal{A} = A$ and $E = \boldsymbol{1} \boldsymbol{1}^T$. Hence,  Theorem~\ref{theorem:main} reduces to
    \[\Psi(A,t)=\left[I - tA \circ / (\boldsymbol{1}\boldsymbol{1}^T - t^2 A \circ A^T) + t^2 \dd(A (A\circ /\boldsymbol{1}\boldsymbol{1}^T - t^2 A \circ A^T))\right]^{-1}\]
thus recovering the known result in the case of a static network \cite{arrigo2024weighted,NQ24}.

\subsection{Updates}\label{sec:updates}

For various reasons, one may wish to update $\mathcal{A}$, the node-level adjacency matrix associated with a time-evolving graph. One natural such update is the addition of a time-stamp, which is in effect the addition of a new block-column, i.e., an update of rank at most $n$, where $n$ is the constant number of nodes in the time-evolving network. This observation can be extremely helpful when backtracking $f$-centralities are to be computed. Indeed, since computing an $f$-centrality amounts to computing $f(\mathcal{A}) \bone$, updating is linked to the problem of updating $f(\mathcal{A})$ after a low-rank update on $\mathcal{A}$; fortunately, this task is well-studied and efficient algorithms exists. See the recent papers \cite{kressner2,kressner1} and the references therein.

In the case of nonbacktracking Katz centrality, however, the theory is not as  well developed because its computation is  based on the algebra of matrices over the ring $R^{N \times N}$ rather than the more familiar $\R^{nN \times nN}$. In Theorem \ref{theorem:update}, we thus analyze the effect of low rank-perturbations (of the special form described above) on the resulting nonbacktracking Katz centrality.

\begin{theorem}\label{theorem:update}
	Define the matrix 
	\[\mathcal{A}_\textit{new}:= \begin{bmatrix}
		\mathcal{A}_\textit{old} & \boldsymbol{A^{[N]}}\\ 0 & A^{[N]}
	\end{bmatrix}\]
as the adjacency matrix associated with the time-evolving graph $\mathcal{G}_\textit{new}$ obtained from $\mathcal{G}$ via the addition of a single final $N$-th time frame, where (here and below) $\boldsymbol{Y}$ denotes a block vector of repeated $Y$ matrices. Further, let $\Psi(\mathcal{A},t)$ be the function described in Theorem~\ref{theorem:main}. Moreover, define $Z_\textit{old} := t\mathcal{A}_\textit{old} \ast (E-t^2 \mathcal{A}_{\textit{old}}^{\ast T} \ast\mathcal{A}_{\textit{old}})^{\ast -1}$, $Z_N :=tA_{[N]} \ast (\boldsymbol{1}\boldsymbol{1}^T-t^2 A_{[N]}^{\ast T} \ast A_{[N]} )^{\ast -1} $ and \begin{align*}
    X &= -Z_\textit{old}  \ast t\mathcal{A}_{\textit{old}}^{\ast T}\ast \boldsymbol{Z_N}- Z_{old} \ast t\boldsymbol{A_{[N]}}^{\ast T} \ast \boldsymbol{Z_N}\\
    & + \textit{dd}^{\ast} \left[t\mathcal{A}_\textit{old}\left(Z_\textit{old}  \ast t\mathcal{A}_{\textit{old}}^{\ast T}\ast \boldsymbol{Z_N}+ Z_{old} \ast t\boldsymbol{A_{[N]}}^{\ast T} \ast \boldsymbol{Z_N} + \boldsymbol{Z_N}\right)\right]  - \boldsymbol{Z_N}  +\textit{dd}^{\ast}\left(t\mathbf{A_{[N]}} Z_N\right).
\end{align*}
 Then,
\[\Psi(\mathcal{A}_\textit{new},t) = \begin{bmatrix}
	\Psi(\mathcal{A}_\textit{old},t) & -\Psi(\mathcal{A}_\textit{old},t)X \boldsymbol{\Psi(A_{[N]},t)}\\ 0 & \Psi(A_{[N]},t)
\end{bmatrix}\]
\end{theorem}

\begin{proof}
	The proof is a matter of expanding the formula of Theorem~\ref{theorem:main}, exploiting where possible the upper-triangular structure of $\mathcal{A}_\textit{new}$.
\end{proof}
\section{Numerical experiments}\label{sec:experiments}

In this section, we test experimentally the increased efficiency of our new methods to compute centrality measures on temporal networks.

\subsection{Timing on randomly generated networks}

In Theorem \ref{thm:fcent} and Corollary \ref{corollary:howtocompute} we showed how to compute, resp., a classical $f$-centrality and nonbacktracking Katz centrality based on the time-evolving adjacency matrix \eqref{eq:teae}
The approaches described in \cite{arrigo2022dynamic,arrigo2024weighted}, involves instead either a larger edge-level adjacency matrix or a weighted Hasmihmoto matrix. The time-evolving adjacency matrix $\mathcal{A}$ \eqref{eq:teae} has size $Nn \times Nn$. For a sparse network, i.e., one with $O(n)$ edges per timerframe, the relevant time-evolving edge-level matrices have also size $O(Nn)$. Nevertheless, \eqref{eq:teae} is typically known, whereas forming the edge-level matrices requires some extra computation. The difference in size is however much more significant for edge-dense graphs, with a number of edges order $O(n^2)$; this situation arises in several applications \cite{noferiniwood2024} including financial mathematics \cite{pozzi2013spread}. In this setting, the edge-level matrices have size $O(Nn^2)$ and the corresponding computational task can quickly become arduous.

\begin{figure}[h]
    
    \begin{subfigure}{.5\textwidth}
    \centering
     \includegraphics[width=\linewidth]{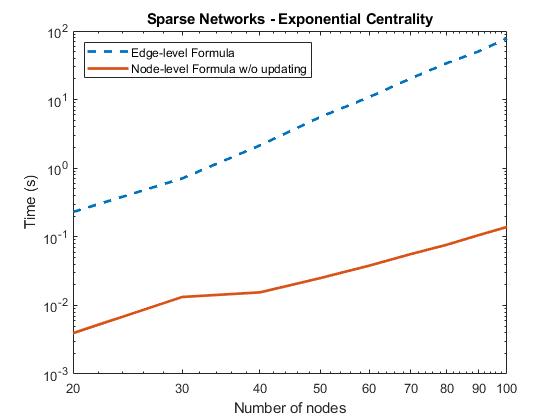}

    \end{subfigure}
    \begin{subfigure}{.5\textwidth}
    \centering
     \includegraphics[width=\linewidth]{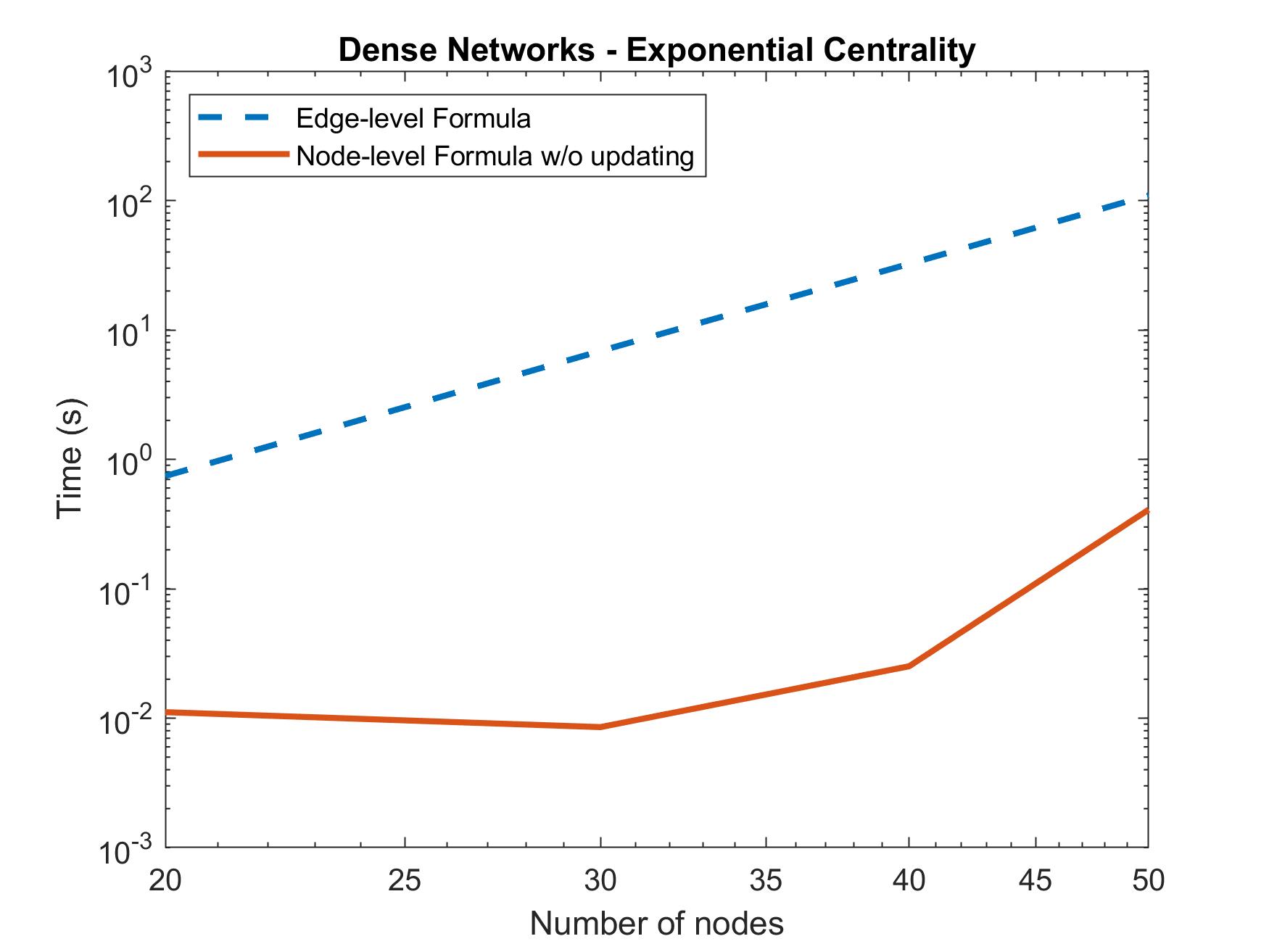}
  
    \end{subfigure}
    \caption{Average computational times for the combinatorially-consistent exponential centrality described in \cite{arrigo2022dynamic,arrigo2024weighted} via the edge-level adjacency matrix (dashed line) versus the time-evolving adjacency matrix (solid line). The average is across 10 different networks with $n   $ nodes and $10$ timeframes, with an expected $3(n-1)$ edges (resp. $3n(n-1)/10$) edges per timeframe for the sparse (resp. dense) networks. }
    \label{fig:f-centrality}
\end{figure}

\begin{figure}[h]
    
    \begin{subfigure}{.5\textwidth}
    \centering
     \includegraphics[width=\linewidth]{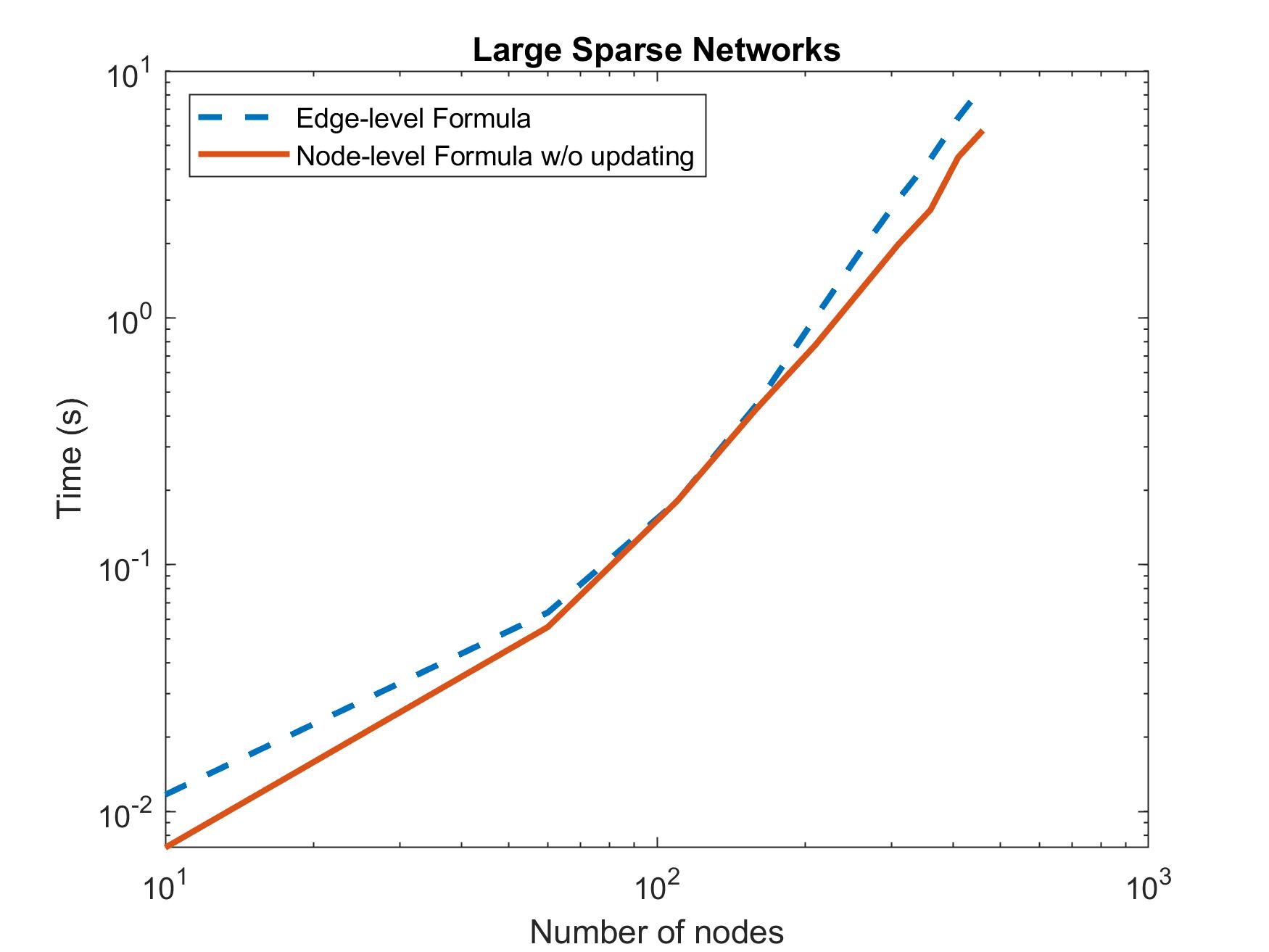}

    \end{subfigure}
    \begin{subfigure}{.5\textwidth}
    \centering
     \includegraphics[width=\linewidth]{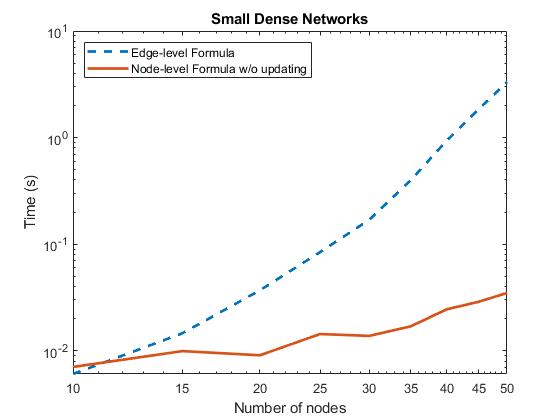}
  
    \end{subfigure}
    \caption{Comparison of the computational time for nonbacktracking Katz centrality using the edge-level approach of \cite{arrigo2022dynamic,arrigo2024weighted} (dashed line) and the new node-level method (solid line), on sparse (left) and dense (right) randomly generated networks with varying number of $n$ nodes and 10 time frames. The expected number of edges is $3(n-1)$ for the sparse networks and $3n(n-1)/10$ for the dense networks.}
    \label{fig:nbtw-centrality}
\end{figure}

\begin{figure}
    \centering
    \includegraphics[width=0.8\linewidth]{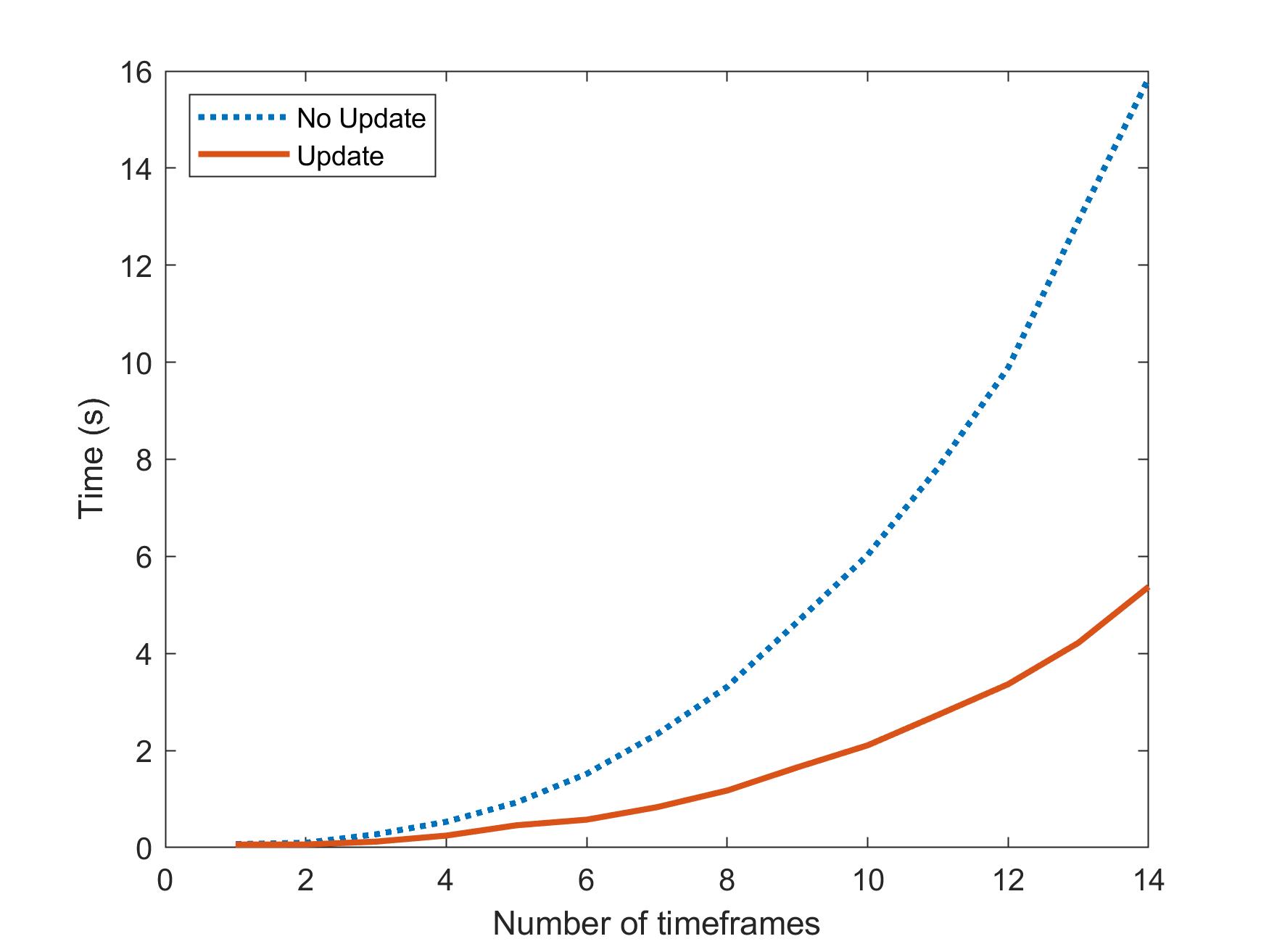}
    \caption{Computation time for a growing return correlation matrix produced by the experiment described in Section \ref{sec:FinanceExperiment} with 480 nodes per time frame. Experimentally, the computational complexities were $O(N^{2.28})$ and $O(N^{2.61})$, resp., for the methods of Corollary \ref{corollary:howtocompute} with and without the update technique described in Theorem~\ref{theorem:update}.} 
    \label{fig:compTimeFrames}
\end{figure}

 This is illustrated in Figure~\ref{fig:f-centrality} and in Figure \ref{fig:nbtw-centrality}, where we examine the computation times, respectively for classical $f$-centralities and nonbacktracking Katz centrality, for sparse and dense randomly generated networks. A further speed-up (not attempted in this experiment) can be achieved by updating the centrality after adding each timeframe; for $f$-centralities this can be done cheaply by leveraging on updating techniques for matrix functions \cite{kressner2,kressner1}, and for nonbacktracking Katz this can be done via Theorem \ref{theorem:update} (see also Figure \ref{fig:compTimeFrames} in the next subsection).

\subsection{Real-world experiment: A network-based method for Portfolio Optimization}\label{sec:FinanceExperiment}
Our development of a new method for computing $f$-centralities and nonbacktracking Katz centrality was motivated by applications where temporal networks have many time frames, and the graphs of each time frame can be very edge-dense. One such scenario arises in financial mathematics, where temporal networks can be used for portfolio optimization in the context of investment management, and these temporal networks typically present precisely the above described features. A full report on the scientific results of this method is beyond the scope of the present paper, and will be presented in a financial mathematics journal. Here, we limit ourselves to a brief summary of the motivating context and of our findings, and we focus instead on the computational aspect, by illustrating how the novel methods of Subection~\ref{sec:NBTDynamic} play a vital role in making the computation of the nonbacktracking Katz centrality tractable within the context of this experiment.

 In \cite{arslan2024portfolio, pozzi2013spread}, a network approach to portfolio construction was investigated. This method involves constructing a static network and investing in certain stocks based on their centralities within this network. In particular, each node represents a stock within a given market. Weighted edges are then inserted between two nodes with a weight based on the correlation of their respective returns over a predetermined time period. This process results in an undirected weighted adjacency matrix associated with the chosen time-period.   It is very reasonable to assume that the relationship governing the correlation of stock returns within a market is time-dependent; however, this information was projected onto a single static network in \cite{arslan2024portfolio, pozzi2013spread} by means of a shifting window. Using the terminology employed in the present paper, this model can be viewed as a time-evolving network in which walks spanning multiple time frames are forbidden. This assumption makes it easier to compute a centrality measure, but does not appear to be justified from the viewpoint of financial mathematics. Thus, it is a natural question to ask whether, by changing the underlying model to a time-evolving network as per Definition~\ref{def:timeevolving}, the financial performance of the method can be improved.

In our experiment, we begin by selecting a time period for analysis and subsequently partition this period into $N$ contiguous intervals. For the $\tau$-th interval, we follow the aforementioned process to produce an adjacency matrix $A^{[\tau]}$ associated with the network $G^{[\tau]}$, where $1\leq \tau \leq N-1$, whereupon we then view the market as a time-evolving network with $\tau$-many time frames, $\mathcal{G} = (G^{[1]}, \dots, G^{[\tau]})$. We then run a simulation where an investor buys $m = 20$ stocks as determined by those with the highest/lowest centrality for a given centrality measure, for the period of the $(\tau+1)$-th interval with portfolio weights obtained solving the Markowitz mean-variance optimization problem, also known as ``Markowitz weights'' \cite{markowitz1952portfolio} and compute the returns yielded by such an investment.  This is then repeated until $\tau = N$, whereupon we sum the returns across the entire period and compare these with the market average for the same period. Each of the $N=O(10)$ time frames consists of 480 stocks. Since each time frame would produce a nearly complete graph, this means that the (dense!) edge-level adjacency matrix of \cite[Definition 5.4]{arrigo2024weighted} has several millions of rows and columns. Computing and storing such a matrix can exceed many computers' memory capabilities. In addition, extending the analysis of \cite{arslan2024portfolio} to temporal networks requires the computation of nonbacktracking Katz for a large number of different combinations of parameters, and hence, the edge-level approach of \cite{arrigo2022dynamic,arrigo2024weighted} is completely unsuitable for this application. Instead, using the node-level approach presented in this paper results in a matrix with just a few thousands rows/columns. This made the computation of nonbacktracking centralities entirely feasible, and allowed us to observe better-than-market returns across a wide range of parameters. Partial results relating to nonbacktracking Katz are presented in brief in Table~\ref{tab:nbtwresults}, and demonstrate that this temporal network approach greatly outperformed the market both in terms of expected return, variance, and Sharpe ratio (computed assuming that the risk-free interest rate is equal to zero).

\begin{table}[]
    \centering
    \begin{tabular}{|c|c|c|}
        \hline  &NBTW-based investment& Market average   \\
          \hline  Expected Return & 20.73 & 12.51\\
          Standard Deviation of Returns & 17.43 & 15.64\\
          Sharpe ratio & 1.148 & 0.326\\
          \hline
    \end{tabular}
    \caption{The best financial performance of the portfolios based on nonbacktracking Katz centrality on temporal networks. }
    \label{tab:nbtwresults}
\end{table}

\bibliographystyle{siam}
        \bibliography{refs}

\begin{thebibliography}{10}

\bibitem{al2021block}
{\sc M.~Al~Mugahwi, O.~D. L.~C. Cabrera, C.~Fenu, L.~Reichel, and
  G.~Rodriguez}, {\em Block matrix models for dynamic networks}, Appl. Math.
  Comput., 402 (2021), p.~126121.

\bibitem{alon2007non}
{\sc N.~Alon, I.~Benjamini, E.~Lubetzky, and S.~Sodin}, {\em Non-backtracking
  random walks mix faster}, Communications in Contemporary Mathematics, 9
  (2007), pp.~585--603.

\bibitem{arrigo2018non}
{\sc F.~Arrigo, P.~Grindrod, D.~J. Higham, and V.~Noferini}, {\em
  Non-backtracking walk centrality for directed networks}, J. Complex Netw., 6
  (2018), pp.~54--78.

\bibitem{arrigo2018exponential}
\leavevmode\vrule height 2pt depth -1.6pt width 23pt, {\em On the exponential
  generating function for non-backtracking walks}, Linear Algebra Appl., 556
  (2018), pp.~381--399.

\bibitem{arrigo2020beyond}
{\sc F.~Arrigo, D.~J. Higham, and V.~Noferini}, {\em Beyond non-backtracking:
  non-cycling network centrality measures}, Proceedings of the Royal Society A,
  476 (2020), p.~20190653.

\bibitem{arrigo2022dynamic}
{\sc F.~Arrigo, D.~J. Higham, V.~Noferini, and R.~Wood}, {\em Dynamic {K}atz
  and related network measures}, Linear Algebra Appl., 655 (2022),
  pp.~159--185.

\bibitem{arrigo2024weighted}
\leavevmode\vrule height 2pt depth -1.6pt width 23pt, {\em Weighted enumeration
  of nonbacktracking walks on weighted graphs}, SIAM J. Matrix Anal. Appl., 45
  (2024), pp.~397--418.

\bibitem{arslan2024portfolio}
{\sc B.~Arslan, V.~Noferini, and S.~Vrontos}, {\em Portfolio management using
  graph centralities: Review and comparison}, arXiv:2404.00187,  (2024).

\bibitem{atiyah}
{\sc M.~F. Atiyah and I.~G. Macdonald}, {\em Commutative algebra},
  Addison-Wesley, Reading, MA (USA), 1969.

\bibitem{kressner2}
{\sc B.~Beckermann, A.~Cortinovis, D.~Kressner, and M.~Schweitzer}, {\em
  Low-rank updates of matrix functions ii: Rational krylov methods}, SIAM J.
  Numer. Anal., 59 (2021), pp.~1325--1347.

\bibitem{kressner1}
{\sc B.~Beckermann, D.~Kressner, and M.~Schweitzer}, {\em Low-rank updates of
  matrix functions}, SIAM J. Matrix Anal. Appl., 39 (2018), pp.~539--565.

\bibitem{bb20}
{\sc M.~Benzi and P.~Boito}, {\em Matrix functions in network analysis}, GAMM
  Mitteilungen, 43 (2020).

\bibitem{benzi2013ranking}
{\sc M.~Benzi, E.~Estrada, and C.~Klymko}, {\em Ranking hubs and authorities
  using matrix functions}, Linear Algebra Appl., 438 (2013), pp.~2447--2474.

\bibitem{boldi2014axioms}
{\sc P.~Boldi and S.~Vigna}, {\em Axioms for centrality}, Internet Mathematics,
  10 (2014), pp.~222--262.

\bibitem{brown}
{\sc W.~C. Brown}, {\em Matrices over commutative rings}, Marcel Dekker, New
  York, NY (USA), 1993.

\bibitem{estrada2012structure}
{\sc E.~Estrada}, {\em The structure of complex networks: theory and
  applications}, American Chemical Society, 2012.

\bibitem{grindrod2018deformed}
{\sc P.~Grindrod, D.~J. Higham, and V.~Noferini}, {\em The deformed graph
  {L}aplacian and its applications to network centrality analysis}, SIAM J.
  Matrix Anal. Appl., 39 (2018), pp.~310--341.

\bibitem{grindrod2011communicability}
{\sc P.~Grindrod, M.~C. Parsons, D.~J. Higham, and E.~Estrada}, {\em
  Communicability across evolving networks}, Phys. Rev. E, 83 (2011),
  p.~046120.

\bibitem{hashimoto1990}
{\sc K.~Hashimoto}, {\em On zeta and l-functions of finite graphs}, Int. J.
  Math., 1 (1990), pp.~381--396.

\bibitem{holme2015modern}
{\sc P.~Holme}, {\em Modern temporal network theory: a colloquium}, The
  European Physical Journal B, 88 (2015), pp.~1--30.

\bibitem{katz1953new}
{\sc L.~Katz}, {\em A new status index derived from sociometric analysis},
  Psychometrika, 18 (1953), pp.~39--43.

\bibitem{kempton2016}
{\sc M.~Kempton}, {\em Non-backtracking random walks and a weighted {I}hara's
  theorem}, Open J. Discrete Math., 6 (2016), pp.~207--226.

\bibitem{laskov2013}
{\sc D.~Laskov}, {\em Diagonalization of matrices over rings}, Journal of
  Algebra, 376 (2013), pp.~123--138.

\bibitem{lin2019non}
{\sc Y.~Lin and Z.~Zhang}, {\em Non-backtracking centrality based random walk
  on networks}, The Computer Journal, 62 (2019), pp.~63--80.

\bibitem{markowitz1952portfolio}
{\sc H.~M. Markowitz}, {\em Portfolio selection}, The Journal of Finance, 7
  (1952), pp.~71--91.

\bibitem{NQ24}
{\sc V.~Noferini and M.~C. Quintana}, {\em Generating functions of
  non-backtracking walks on weighted digraphs: {R}adius of convergence and
  {I}hara's theorem}, Linear Algebra Appl., 699 (2024), pp.~72--106.

\bibitem{noferiniwood2024}
{\sc V.~Noferini and R.~Wood}, {\em {Efficient computation of Katz centrality
  for very dense networks via negative parameter Katz}}, J. Complex Netw., 12
  (2024), p.~cnae036.

\bibitem{pozzi2013spread}
{\sc F.~Pozzi, T.~Di~Matteo, and T.~Aste}, {\em Spread of risk across financial
  markets: better to invest in the peripheries}, Sci. Rep., 3 (2013), p.~1665.

\bibitem{strogatz2001exploring}
{\sc S.~H. Strogatz}, {\em Exploring complex networks}, nature, 410 (2001),
  pp.~268--276.

\bibitem{timar2021approximating}
{\sc G.~Tim{\'a}r, R.~Da~Costa, S.~Dorogovtsev, and J.~Mendes}, {\em
  Approximating nonbacktracking centrality and localization phenomena in large
  networks}, Phys. Rev. E, 104 (2021), p.~054306.

\end{thebibliography}

        \appendix
        
        \section{A spectral theory for matrices of vectors}\label{sec:matvec}

In this appendix, we develop a spectral theory of matrices whose elements lie in a ring that is a finite-dimensional vector space over a field, where (having fixed a basis) addition and multiplication are defined entrywise. This algebraic detour is helpful for the development of a fast method to compute nonbacktracking centrality on temporal networks. In particular, Subsection \ref{sec:NBTDynamic} relies on the special case of block matrices whose block elements lie in the commutative ring of square real matrices enhanced with the usual matrix addition and the so-called Schur (elementwise) product.

\subsection{Determinant, matrix multiplication, identity, matrix inverse}\label{sec:dee}

Let $R$ be a nontrivial commutative ring with unity. Then $R^{m \times m}$ is the (non-commutative) ring of $m \times m$ square matrices over $R$, where matrix addition is defined elementwise and matrix multiplication is built in the usual way from the two basic operations in $R$. The determinant is a function from $R^{m \times m}$ to $R$, defined in the usual way and denoted by $\det_R$, and it still satisfies most of the familiar properties that it possesses over a field: see \cite{brown} for more details. However, the notions of invertible and nonsingular matrix no longer necessarily coincide.
\begin{definition}
     A square matrix $M \in R^{m \times m}$ is \emph{invertible over $R$} if there exists $N\in R^{m \times m}$, called the inverse of $M$, such that $M N = N M = E$; and it is \emph{singular} if its kernel is notrivial, i.e., there exist $0 \neq v \in R^m$ such that $Mv=0$.
\end{definition}
    One can prove \cite[Corollary 2.21]{brown} that $M$ is invertible over $R$ if and only if $\det M \in R^\times$ is a unit of $R$; and that $M$ is singular if and only if $\det M$ is a zero divisor of $R$. Note that in general it is possible that a matrix is neither singular nor invertible over $R$, unless every element of $R$ is either a unit or a zero divisor.

Henceforward, we specialize to the commutative ring $R=\{ \mathbb{F}^{p}, +, \circ \}$, that is, $p$-vectors over the algebraically closed field $\F$ enhanced with addition and \emph{elementwise} multiplication. The additive identity of $R$ is the zero vector, whereas the multiplicative identity is the vector of all ones, which we denote by $\bone$. While the theory developed in this section is valid for any $p$ and any closed field $\F$, for the purposes of this paper we are particularly interested in the special case of square matrices with complex entries. To this end, we set\footnote{As a matter of fact, the matrices of our interest in this paper are all real; but as real matrices may have nonreal eigenvalues, it is necessary for our developments to embed $\R \subseteq \C$.} $\F=\C$ and $p=n^2$, and further make the identification  $\F^{n^2} \cong \F^{n \times n}$; that is, we take as our commutative ring $R=\{ \C^{n \times n},+,\circ\}$. In this special case, the additive identity is the zero matrix and the multiplicative identity is the rank-one matrix of all ones, $\bone \bone^T$. $R=\{ \mathbb{F}^{p}, +, \circ \}$ is also a finite-dimensional vector space over $\F$, and hence an Artin ring \cite[Ch. 8]{atiyah}. This in turn implies that (i)  there are only finitely many (more precisely, $2^{p}$) ideals in $R$; (ii) every element of $R$ is either a unit (if it has no zero entries) or a zero divisor (if it has at least one zero entry).

Let us consider the set of $m \times m$ square matrices over $R$; matrix addition and multiplication can be defined as usual, by leveraging scalar addition and multiplication in $R$. We thus obtain a noncommutative ring, whose additive identity is the zero block matrix and whose multiplicative identity is $E=\bigoplus_{i=1}^m \bone$  (for $R=\{\F^p,+,\circ\}$), or $E=\bigoplus_{i=1}^m \bone \bone^T$ (when $R=\{\F^{n \times n},+,\circ\}$).
It is worth noting that, in the latter case, an element of $R^{m \times m}$ can also represent an element of $\F^{mn \times mn}$, but the definitions of matrix multiplication differ accordingly. To avoid confusion, we shall denote matrix multiplication in $\F^{mn \times mn}$ by the usual juxtaposition of symbols, and matrix multiplication in $R^{m \times m}$ by the symbol $\ast$: hence, if $A,B \in R^{m \times m}$, we have 
\[ (A \ast B)_{rs} = \sum_{k=1}^m A_{rk} \circ B_{ks} \in R \]
where $\circ$ denotes the entrywise product of matrices, sometimes also called the Schur or Hadamard product. It is clear that in general $A \ast B \neq AB$. Furthermore given $\lambda \in R$ and either $A \in R^{m \times m}$ or $A \in R^m$, we write $\lambda \circ A$ to mean element-wise multiplication by the constant $\lambda$, i.e., $[\lambda \circ A]_{rs} = \lambda \circ A_{rs}$ or $[\lambda \circ A]_{r} = \lambda \circ A_{r}$ respectively.

 If $M$ is invertible over $R$, we denote its inverse (over $R$) by $M^{\ast -1}$ so that $M \ast M^{\ast -1} = M^{\ast -1} \ast M=E$. Finally, positive powers of a matrix $M$ are defined as $M^{\ast k}=\underbrace{M \ast M \ast \dots \ast M}_{k \ \mathrm{times}}$; moreover, we define $M^{\ast 0}=E$, and for matrices invertible over $R$ negative powers are defined as $M^{\ast -k}:=(M^{\ast -1})^{\ast k}$.
 \subsection{Eigenvalues and eigenvectors}\label{sec:ee}
 Clearly, we have the isomorphisms \[ R \cong \underbrace{\F \times \F \times \dots \times \F}_{p \ \mathrm{times}} \Rightarrow R^{m \times m} \cong \underbrace{\F^{m \times m} \times \F^{m \times m}  \times \dots \times \F^{m \times m} }_{p \ \mathrm{times}}.\] In particular, given $A \in R^{m \times m}$ and $1 \leq i \leq p$, we can define $[A]_{i} \in \F^{m \times m}$ as $([A]_{i})_{rs} = (A_{rs})_{i}$; in the important special case where $R$ is the ring of $n \times n$ complex matrices with elementwise addition and multiplication, this definition constructs $n^2$ matrices $[A]_{ij} \in \C^{m \times m}$ as $([A]_{ij})_{rs} = (A_{rs})_{ij}$. Having set up this notation, the above observation immediately implies some interesting properties.
 \begin{corollary}\label{cor1}
     When $R=\{\F^p,+,\circ\}$ then, for all $A,B \in R^{m \times m}$,
\[ ((A \ast B)_{rs})_{i} = ([A]_{i}[B]_{i})_{rs}.  \]
for all $i=1,\dots,p$. In particular, when $R=\{\C^{n \times n},+,\circ\}$ and for all $1 \leq i,j \leq n$ we have
\[ ((A \ast B)_{rs})_{ij} = ([A]_{ij}[B]_{ij})_{rs}.  \]
 \end{corollary}
Note that on the left-hand side of the equations in Corollary \ref{cor1} we have matrix multiplication over the ring $R$, but on the right-hand side we have the \emph{usual} matrix multiplication over the field $\F$.
\begin{corollary}\label{cor2}
    For all $A \in R^{m \times m}$ and all $i=1,\dots, p$; $(\det_R A)_{i} = \det [A]_{i}$. In particular, when $R=\{\C^{n \times n},+,\circ\}$ we have $(\det_R A)_{ij} = \det [A]_{ij}$ for all $1 \leq i,j \leq n$.
\end{corollary}
In the equations in Corollary \ref{cor2}, $\det$ denotes the usual determinant defined for square matrices over a field, while $\det_R$ is the determinant over $R$.

Remarkably, even though the underlying ring is not an integral domain, one can still build a sensible theory of eigenvalues for matrices in $R^{m \times m}$. To this goal, we define the eigenvalues of a square matrix over $R$ via the notion of a characteristic polynomial: having fixed $M \in R^{m \times m}$, this is defined as the function 
\[p_M : R \rightarrow R, \qquad z \in R \mapsto p_M(z)=\mathrm{det}_{R} (z \circ E - M).\]
\begin{definition}\label{defeig}
    An element $\lambda \in R$ is an eigenvalue of $M \in R^{m \times m}$ if it is a root the characteristic polynomial of $M$, i.e., if $p_M(\lambda)=\det_R(\lambda \circ E-M)=0$.
\end{definition}
The number of eigenvalues of a square matrix over $R$, per Definition \ref{defeig}, is finite albeit quite large. Sensible notions of eigenvectors and  diagonalization \cite{laskov2013} also emerge.

\begin{theorem}\label{thmdiag}
Let $R=\{ \F^{n \times n},+,\circ \}$.
    \begin{enumerate}
    \item Every matrix $M \in R^{m \times m}$ has $m^{n^2}$ eigenvalues according to Definition \ref{defeig}, counted with algebraic multiplicities;
        \item If $\lambda$ is an eigenvalue of $M$ according to Definition \ref{defeig}, then there is a nonzero vector $V \in R^m$, called an eigenvector, such that $M \ast V=\lambda \circ V$ (but, generally, the reverse implication is false);
        \item Suppose further than, for all $i,j=1,\dots,n$, the matrices $[M]_{ij}$ are diagonalizable. Then, there exists a choice of $m$ eigenvalues of $M$, say, $\lambda_1,\dots,\lambda_m$, that are associated with a linearly independent set of eigenvectors $V_1,\dots,V_m$. In particular, $M$ is diagonalizable over $R$ in the sense that $M\ast V=V\ast \Lambda$ with
        $ V=\begin{bmatrix}
            V_1 & \dots & V_m
        \end{bmatrix} \in R^{m \times m}, \qquad  
 \Lambda = \bigoplus_{i=1}^m \lambda_i \in R^{m \times m}$, and $V$ is invertible over $R$ so that $M=V \ast \Lambda \ast V^{\ast -1}$.
    \end{enumerate}
\end{theorem}
\begin{proof}
\begin{enumerate}
\item    By Corollary \ref{cor2},
    \[  \mathrm{det}_R (\lambda \circ E - M) = 0 \Leftrightarrow \det(\lambda_{i,j} I - [M]_{i,j})=0 \ \forall \ i,j .\]
Hence, the $(i,j)$ element of $\lambda \in R=\F^{n \times n}$ must be an eigenvalue (in the classical sense) of $[M]_{i,j} \in \F^{m \times m}$. There are $m$ possible choices (counting multiplicities) for each element, and hence there are $m^{n^2}$ eigenvalues of $M$.
    \item Let $\lambda$ be an eigenvalue of $M$; then for all pairs $(i,j)$ it holds $[M]_{i,j} v_{i,j} = \lambda_{ij} v_{i,j}$ for some nonzero vector $v_{i,j} \in \F^m$. Define $V \in R^m$ so that $(V_k)_{i,j}=(v_{i,j})_k$: then, $V\neq0$ and $M \ast V = \lambda \circ V$ by Corollary \ref{cor1}. To see that the reverse implication is generally false, observe that $V \neq 0$ and $M \ast V = \lambda \circ V$ implies that $\lambda \circ E - M$ is singular, which is a weaker requirement than its determinant being $0$. In particular, if $0 \neq p_M(\lambda)$ is singular, then $\lambda$ is associated with a nonzero $V \in R^m$ such that $M \ast V=\lambda \circ V$, but is not an eigenvalue in the sense of Definition \ref{defeig}. Note, in particular, that if $\F=\C$ and $n>1$ then generally there are uncountably many values of $\lambda$ that satisfy $M \ast V = \lambda \circ V$ for some nonzero $V$, but only finitely many that satisfy $p_M(\lambda)=0$.
    \item By assumption, for all pairs $(i,j)$ there exist a diagonal matrix $D_{ij}$ and an invertible (over $\F$) matrix $W_{ij}$ such that $[M]_{ij} = W_{ij} D_{ij} W_{ij}^{-1}$. Now, define $V,\Lambda \in R^{m \times m}$ by the relations $[V]_{ij}=W_{ij}$, $[\Lambda]_{ij}=D_{ij}$. Then: (a) $\Lambda$ is diagonal over $R$ (that is, block diagonal over $\F$ with $n \times n$ blocks); (b) $V$ is invertible over $R$ by Corollary \ref{cor2} because $(\det_R V)_{ij} = \det W_{ij} \neq 0$; (c) $M \ast V = V \ast \Lambda$ follows by Corollary \ref{cor1} (d) Denoting $\lambda_i := \Lambda_{ii}$ we have that $\det_R(\lambda_i \circ E-M)=\det_R(\lambda_i \circ E-\Lambda)=0$. Moreover, $M \ast V_i = \lambda_i \circ V_i$ where $0 \neq V_i \in R^m$ is the $i$-th column of $V$.
\end{enumerate}
\end{proof} 
    Generally, the $m$-uple of eigenvalues $\lambda_k$ in Theorem \ref{thmdiag}, item 3, cannot be arbitrarily selected from the $m^{n^2}$ eigenvalues of $M$.
    In the generic case where every $[M]_{ij}$ has $m$ distinct eigenvalues, the number of the possible $m$-uples of eigenvalues that diagonalize $M$ is (up to permutations) $(m!)^{n^2-1}$, which is considerably lower than $\begin{pmatrix}
        m^{n^2}\\
        m
    \end{pmatrix}$, i.e., the number of distinct $m$-uples from a set with $m^{n^2}$ elements.

If an $m$-uple is selected arbitrarily rather than following the construction in the proof of item 3, the statement of item 3 may fail. For a counterexample, take $\F=\C$ and $M=\begin{bmatrix}
    A&0\\
    0&-A
\end{bmatrix}$ with
$A=\begin{bmatrix}
    1&2\\
    3&4
\end{bmatrix}.$ Then,
$ [M]_{11}=\begin{bmatrix}
    1 &0 \\
     0& -1
\end{bmatrix}, \ [M]_{12}=\begin{bmatrix}
    2 &0 \\
     0& -2
\end{bmatrix}, \ [M]_{21}=\begin{bmatrix}
    3 &0 \\
     0& -3
\end{bmatrix}, \ [M]_{22}=\begin{bmatrix}
    4 &0 \\
     0& -4
\end{bmatrix}.$ Consider the eigenpairs $(\lambda_1,V_1)$ and $(\lambda_2,V_2)$ with
$\lambda_1 = \begin{bmatrix}
    1 & 2\\
    3 & 4
\end{bmatrix}, 0 \neq V_1=\begin{bmatrix}
    X\\
    0
\end{bmatrix}; \qquad \lambda_2 = \begin{bmatrix}
    1 & 2\\
    3 & -4
\end{bmatrix}, 0 \neq V_2 = \begin{bmatrix}
    \begin{bmatrix}
        x & y\\
        z &0
    \end{bmatrix}\\
    \begin{bmatrix}
        0 & 0\\
        0 & w
    \end{bmatrix}
\end{bmatrix}.$ Then,  $\mathrm{det}_R \begin{bmatrix}
    V_1 & V_2
\end{bmatrix} = X\circ \begin{bmatrix}
    0 & 0\\
    0 & w
\end{bmatrix}  \not \in R^\times$ implying that $V_1$ and $V_2$ are linearly dependent regardless of the choice of $0 \neq X \in R$ and of $x,y,z,w \in \C$ not all zero.

\begin{remark}\label{rem:jordan}
    When some of the matrices $[M]_{ij}$ are not diagonalizable, it is possible to modify item 3 in Theorem \ref{thmdiag} by removing that assumption. Rather than forming the matrix $V$ with eigenvectors of $[M]_{ij}$, one can start from (still linearly independent) sets of generalized eigenvectors, i.e., Jordan  chains. This still allows us to write $M \ast V=V\ast \Lambda$, but $\Lambda$ is no longer (block) diagonal. It is instead in ``generalized Jordan form": the superdiagonal blocks are allowed to have elements in $\{0,1\}$. For example, setting $n=2$, then the matrices $\Lambda_1,\Lambda_2 \in R^{3 \times 3} $ where 
\[  \Lambda_1 = \begin{bmatrix}
    \begin{bmatrix}
        1&2\\
        3&4
    \end{bmatrix} & \begin{bmatrix}
        1&1\\
        0&0
    \end{bmatrix} & \begin{bmatrix}
        0&0\\
        0&0
    \end{bmatrix}\\
    \begin{bmatrix}
        0&0\\
        0&0
    \end{bmatrix}&\begin{bmatrix}
        1&2\\
        3&4
    \end{bmatrix}&\begin{bmatrix}
        0&0\\
        0&0
    \end{bmatrix}\\
    \begin{bmatrix}
        0&0\\
        0&0
    \end{bmatrix}&\begin{bmatrix}
        0&0\\
        0&0
    \end{bmatrix}&\begin{bmatrix}
        1&2\\
        3&4
    \end{bmatrix}
\end{bmatrix}, \Lambda_2 = \begin{bmatrix}
    \begin{bmatrix}
        1&2\\
        3&4
    \end{bmatrix} & \begin{bmatrix}
        1&1\\
        1&1
    \end{bmatrix} & \begin{bmatrix}
        0&0\\
        0&0
    \end{bmatrix}\\
    \begin{bmatrix}
        0&0\\
        0&0
    \end{bmatrix}&\begin{bmatrix}
        1&2\\
        3&4
    \end{bmatrix}&\begin{bmatrix}
        1&0\\
        0&1
    \end{bmatrix}\\
    \begin{bmatrix}
        0&0\\
        0&0
    \end{bmatrix}&\begin{bmatrix}
        0&0\\
        0&0
    \end{bmatrix}&\begin{bmatrix}
        1&2\\
        3&4
    \end{bmatrix}
\end{bmatrix}  \]
    are both in ``generalized Jordan form".
    
\end{remark}

While we have stated  Theorem \ref{thmdiag} for the special case $R=\{ \F^{n \times n},+,\circ \}$, it is clear that it holds more generally for $R=\{ \F^p,+,\circ \}$; we state this as Theorem \ref{thmdiag2}.

\begin{theorem}\label{thmdiag2}
    Let $R=\{ \F^{p},+,\circ \}$.
    \begin{enumerate}
    \item Every matrix $M \in R^{m \times m}$ has $m^{p}$ eigenvalues according to Definition \ref{defeig}, counted with algebraic multiplicities;
        \item If $\lambda$ is an eigenvalue of $M$ according to Definition \ref{defeig}, then there is a nonzero vector $V \in R^m$, called an eigenvector, such that $M \ast V=\lambda \circ V$ (but, generally, the reverse implication is false);
        \item Suppose further than, for all $i=1,\dots,p$, the matrices $[M]_{i}$ are diagonalizable. Then, there exists of choice of $m$ eigenvalues of $M$, say, $\lambda_1,\dots,\lambda_m$, that are associated with a linearly independent set of eigenvectors $V_1,\dots,V_m$. In particular, $M$ is diagonalizable over $R$ in the sense that $M\ast V=V\ast \Lambda$ with
         $V=\begin{bmatrix}
            V_1 & \dots & V_m
        \end{bmatrix} \in R^{m \times m}, 
 \Lambda = \bigoplus_{i=1}^m \lambda_i \in R^{m \times m}$, and $V$ is invertible over $R$ so that $M=V \ast \Lambda \ast V^{\ast -1}$.
    \end{enumerate}
\end{theorem}
\subsection{Convergence of infinite series}

Let us now specialize to the case $\F=\C$, and endow $R^m \cong \C^{m n^2}$ (or, more generally, $R^m \cong \C^{pm}$) with an arbitrary norm $\| \cdot \|$. This in turn induces an operator norm on $R^{m \times m}$ via
\[  \| A \|_{op} = \sup_{x:\|x\|=1} \|A \ast x \|. \] It is easy to verify that $\| A \ast B \|_{op} \leq \|A \|_{op} \|B\|_{op}$. All these (vector and matrix) norms are equivalent since they ultimately are norms on a finite-dimensional vector space over $\C$. Thus, we henceforth assume that one arbitrary choice of norm on $R^m$ and  a corresponding operator norm on $R^{m \times n}$ has been made.

Given $M \in R^{m \times m}$, a scalar $z \in \C$, and a sequence $(c_k)_k \subset \C$, we are interested in studying the convergence of the infinite series $\sum_{k=0}^\infty c_k z^k M^{\ast k}$; convergence here means that there is a limit $L \in R^{m \times m}$ such that for all $\epsilon > 0$ there is $N \in \N$ such that
\[ n \geq N \Rightarrow \left\| L - \sum_{k=0}^n c_k z^k M^{\ast k}\right\|_{op} < \epsilon.    \]
Similarly, the sequence is divergent if for all $\delta \in \R$ there is $N \in \N$ such that
\[ n \geq N \Rightarrow \left\| L - \sum_{k=0}^n c_k z^k M^{\ast k}\right\|_{op} > \delta.    \]

\begin{theorem}\label{thm:68}
    Let $R=\{ \C^{n \times n}, +, \circ \}$, and let $M \in R^{m \times m}$. If the scalar power series $\sum_{k=0}^\infty c_k z^k$ has radius of convergence $r$ and we define the generalized spectral radius as
    \[ \rho_R(M) := \max_{i,j} \rho([M]_{ij}), \]
    where $\rho(X)$ is the usual spectral radius of the square matrix $X$ over $\C$, then the series 
    $\sum_{k=0}^\infty c_k z^k M^{\ast k}$ is convergent for $|z| < \frac{r}{\rho_R(M)}$ and it is divergent for $|z| > \frac{r}{\rho_R(M)}.$

The same result holds for $R=\{\C^p,+,\circ\}$, but in this case the generalized spectral radius is
   \[ \rho_R(M) = \max_{i} \rho([M]_{i}).\]
\end{theorem}
\begin{proof}
    We prove the result for $R=\{\C^{n \times n},+,\circ\}$; the more general case is analogous. Suppose first that the matrices $[M]_{ij}$ are diagonalizable. By Theorem \ref{thmdiag}, $ M^{\ast k} = V \ast \Lambda^{\ast k} \ast V^{\ast -1}$. It follows that the radius of convergence of the sequence of partial sums $\sum_{k=0}^n c_k z^k M^{\ast k}$ is equal to the minimum of the radii of convergence of the $m n^2$ scalar sequences $\sum_{k=0}^n c_k z^k (\lambda_{ij,\ell})^k$ where $\lambda_{ij,\ell}$ ($\ell=1,\dots,m$) are the eigenvalues of $[M]_{ij}$. The statement follows by observing that $|\lambda_{ij,\ell}| \leq \rho([M]_{ij})$ and equality is achieved by at least one value of $\ell$.

More generally, some of the matrices $[M]_{ij}$ might have non-diagonal Jordan form. Recall that, if $J$ is a Jordan block with eigenvalue $\lambda$, then the radius of convergence of the sequence of partial sums $\sum_{k=0}^n c_k z^k J^k$ is equal to $r/\lambda$.
Using Remark \ref{rem:jordan}, we conclude that the statement holds also in this case because, by the isomorphism described in Section \ref{sec:ee}, the scalar entries of $\Lambda^k$ are the elements of $J_{ij\ell}^k$ where $J_{ij\ell}$ is the $\ell$-th Jordan block in an (arbitrary, but fixed) Jordan canonical form of $[M]_{ij}$. 
\end{proof}
\end{document}